%% file: SDE_jumps_6_9_2006.tex
\documentclass[10pt]{amsart}
\usepackage{amssymb,amsmath}
\usepackage{graphics}
\usepackage[hypertex]{hyperref}
\usepackage{hyperref}
\newtheorem{thm}{Theorem}[section]
\newtheorem{lemma}[thm]{Lemma}

\newtheorem{defin}[thm]{Definition}
\newtheorem{remark}[thm]{Remark}
\numberwithin{equation}{section}
\newcommand\e{\mathbf{e}}
\newcommand{\ssy}{\scriptscriptstyle}
\newcommand{\bara}{\overline{a}}
\newcommand{\barb}{\overline{b}}
\newcommand{\bard}{\overline{d}}
\newcommand{\barX}{\overline{X}}
\newcommand{\barY}{\overline{Y}}
\newcommand{\barN}{\overline{N}}
\newcommand{\Psistar}{{\mathcal P}^{n+1}}
\newcommand{\Rstar}{{\widehat{\mathcal P}}^n}
\newcommand{\ttn}{{m}}
\newcommand{\ttt}{\tilde{t}}
\newcommand{\dtt}{\Delta\tilde{t}}
\newcommand{\sigmatt}{\mathcal{G}} 

\newcommand{\Z}{{\bf Z}}
\newcommand{\cF}{\mathcal{F}}
\newcommand{\cB}{\mathcal{B}}
\newcommand{\cJ}{\mathcal{J}}
\newcommand{\OM}{\mathcal{O}}

\newcommand{\brr}{\overline}
\newcommand{\qqq}{\partial}
\newcommand{\tol}{\mathrm{TOL}}
\newcommand{\ETS}{\mathtt{E}_{TS}}
\newcommand{\rset}{\mathbb{R}}
\newcommand{\nset}{\mathbb{N}}
\newcommand{\Ito}{\mbox{It{\^o}} }

\newcommand{\supdt}{\Delta t_{sup}}
\newcommand{\var}{\text{\rm Var}}
\newcommand{\ErrS}{\mathtt{E}_S}
\newcommand{\mF}{\mu} 
\newcommand{\ma}{\mathcal{A}}
\newcommand{\mc}{\mathcal{C}}
\newcommand{\me}{\mathcal{E}}

\newcommand{\nA}{\Phi} 
\newcommand{\normal}{\mathcal{N}}
\newcommand{\pa}{\partial}
\newcommand{\ra}{\rightarrow}

\newcommand{\f}{\frac}
\newcommand{\ken}{\ \ }
\newcommand{\kenn}{\hskip0.1truecm}

\newcommand{\half}{\frac{1}{2}}
\font\aut=cmr8
\def\errorc{\mathcal{E}_{\ssy C}}
\def\errort{\mathcal{E}_{\ssy T}}
\def\errors{\mathcal{E}_{\ssy S}}

\begin{document}
\title[]{Adaptive Weak Approximation of Diffusions with Jumps}
\author[]
{E. Mordecki$^{\dag}$, A. Szepessy$^{\ddag}$,
R. Tempone$^{\S}$ and G. E. Zouraris$^{\P}$}
\thanks{%
$^{\dag}$Centro de Matem\'atica, Facultad de Ciencias,
Universidad de la Rep\'ublica, Igu\'a 4225 C.P. 11400, Montevideo,
Uruguay (mordecki@cmat.edu.uy).}
\thanks{%
$^{\ddag}$Matematiska Institutionen,
Kungl. Tekniska H{\"o}gskolan,
S--100 44 Stockholm, Sweden
(szepessy@nada.kth.se).}
\thanks{%
$^{\S}$Department of Mathematics and School of Computational
Science and Information Technology, 400 Dirac Science Library,
Florida State University, Tallahassee, FL 32306--4120
(rtempone@scs.fsu.edu).}
\thanks{%
$^{\P}$Department of Mathematics, University of Crete,
GR--714 09 Heraklion, Crete, Greece and Institute of
Applied and Computational Mathematics, FO.R.T.H.,
GR--711 10 Heraklion, Crete, Greece (zouraris@math.uoc.gr).}
%
%
\begin{abstract}
This work develops Monte Carlo Euler adaptive time
stepping methods for the weak
approximation problem of jump diffusion driven stochastic
differential equations.
The main result is the derivation of a new expansion
for the computational error, with computable leading
order term in a posteriori form, based on stochastic
flows and discrete dual backward problems which
extends the results in \cite{STZ}. These expansions
lead to efficient and accurate computation of error
estimates.
Adaptive algorithms for either
stochastic time steps or quasi-deterministic
time steps are described.
Numerical examples show the performance of
the proposed error approximation and of the
described adaptive time-stepping methods.
\end{abstract}
\maketitle
\par\noindent\medskip\par\noindent
\centerline{\today}
%
%
%
\section{Introduction}\label{sec:intro}
This work develops adaptive methods and proves a posteriori error expansions,
with computable leading order term, for weak approximation of jump-diffusions
driven stochastic differential equations.
\subsection{Problem's setting}

Consider $X=\{X(t)=(X^1(t),\dots,X^d(t))\colon t\in[0,T]\}$,
a $d$-dimensional stochastic process
that is the solution of the stochastic differential equation
%
%
%
%
\begin{equation}\label{eq:1.1}
\begin{split}
X(t)=X(0)+&\int_0^ta(\,s,X(s^{-})\,)\,ds
+\sum_{\ell=1}^{\ell_0}\int_0^t b^{\ell}(\,s,X(s^{-})\,)
\,dW^{\ell}(s)\\
+&\int_0^t\int_{\Z}c(\,s,X(s^{-}),z\,)\,p(ds,dz),
\end{split}
\end{equation}
on a time interval $[0,T]$ (see III.2c in \cite{js}, or,
\cite{IW}). The randomness in the equation is generated by (i) an
$\cF_0-$measurable $d$-dimensional random variable $X(0)$,
(ii) a standard $\ell_0$-dimensional Wiener process
$W=\{W(t)=(W^1(t),\dots,W^{\ell_0}(t))\colon t\in[0,T]\}$,
i.e. its coordinates are independent standard real
valued Wiener processes,
and (iii) a Poisson random measure $p$ on $[0,T]\times \Z$, where
$\Z\equiv\rset^{\ell_1}\backslash\{0\}$,
with deterministic time dependent intensity measure
$q(dt,dz)=\lambda(t)dt\otimes \mF(t,dz)$. 
Here $\lambda\colon[0,T]\to [0,\lambda_{\rm max}]$
is the time intensity of jumps with $\lambda_{\rm max}\in(0,+\infty)$
and, for each $t\geq 0$ fixed, $\mF(t,dz)$ is a probability measure on $\Z$.
All processes are defined on a stochastic basis
$\cB = (\Omega,\cF_T,\{\cF_t\}_{0\le t\le T}, P)$,
and, as usual in this framework, we assume that the Wiener process
and the Poisson random measure are independent.
The coefficients
$a\colon[0,T]\times\rset^d\rightarrow\rset^d$,
$b\colon[0,T]\times\rset^d\rightarrow\rset^{d \times \ell_0}$
and $c\colon[0,T]\times\rset^d\times\Z\to\rset^d$
are assumed to be Borel functions, satisfying regularity
conditions
defined in Lemma~\ref{lem:2.1}.
\par
For a given scalar function $g:\rset^d\rightarrow\rset$,
the goal of our work is to construct approximations to the expected
value $E\big[\,g\big(X(T))\big]$ by a Monte Carlo Euler method
(cf., e.g.,  \cite{KP}, \cite{Mil}).
In what respects the dimension, although our developments
are valid for all $d\in{\mathbb N}$, we have in mind relatively
high values of $d$, taking into account the curse of dimensionality
present in concurrent deterministic methods that solve partial integro-differential
equations.
The computational approach presented is common in computing option prices
in mathematical finance and in simulating stochastic dynamics
(cf., e.g.,  \cite{CT_book}, \cite{KP}, \cite{KS}, \cite{O}).
%
%
%
%
%
%
\begin{remark}[Construction of the integral with respect to the Poisson
random measure]\label{Remarkk_1}
\mbox{}
\par\noindent
Consider a sequence $\e_1,\e_2,\dots$ of independent random variables
with common exponential distribution with parameter 1. Define
\begin{equation}\label{eq:lambda_def}
\Lambda(t)=\int_0^t\lambda(s)\,ds,\quad t\in[0,T].
\end{equation}
The number of jumps of the random Poisson measure $p(dt,dz)$ in an interval $[0,t]$
is determined as
\begin{equation*}
\widehat{N}(t)=\max\Big\{k\colon \sum_{j=1}^k\e_j\leq \Lambda(t)\Big\},
\end{equation*}
and the total number of jumps in $[0,T]$
is denoted by $\widehat{N}\equiv\widehat{N}(T)$.
The jump times of the Poisson measure can be defined
by $\tau_0=0$, and $\tau_{k}=\Lambda^{-1}(\e_1+\cdots+\e_k)$ for $k=1,\dots,\widehat{N}$,
and can be computed recursively by
\begin{equation}\label{eq:tau_def}
\e_k=\int_{\tau_{k-1}}^{\tau_{k}}\lambda(s)\,dt,\quad\text{for}\quad k=1,\dots,\widehat{N}.\\
\end{equation}
Once the jump times are computed, we proceed to sample the marks
$\{Z_{k}\}$, that, conditionally on the values of
the jumps times, are independent random variables distributed respectively according to
$\{\mF(\tau_{k},dz)\}$.
The random measure with intensity $q(dt,dz)=\lambda(t)dt\otimes\mu(t,dz)$
can then be constructed as
$$
p(ds,dz)=\sum_{k=1}^{\widehat{N}}\delta_{(\tau_{k}, Z_{k})}(ds,dz),
$$
and, consequently, the stochastic integral with respect to
the Poisson random measure, i.e. the last term in \eqref{eq:1.1},
can be computed as
\begin{equation*}
\int_0^t\int_{\Z}c(\,s,X(s^{-}),z\,)\,p(ds,dz)
=\sum_{k=1}^{\widehat{N}(t)}c\big(\tau_{k},X(\tau_{k}^{-}),Z_{k}\big),\quad t\in[0,T].
\end{equation*}
\end{remark}
%
%
\begin{remark}\label{rem:EulerI}
Since the time intensity $\lambda(t)$ is deterministic,
the jump times $\{\tau_{k}\}$ can
be directly obtained before solving for the process $X$;
it is enough then to find the function $\Lambda$ by accurately
integrating \eqref{eq:lambda_def} and then successively finding
$\{\tau_{k}\}$ by solving \eqref{eq:tau_def}.
Furthermore, if $\lambda$ is just a constant then we simply
sample $\tau_{k}=\e_1+\cdots+\e_k$ from a sequence
$\{\e_k\}$
of independent random variables
with common exponential distribution with parameter $\lambda$.
\end{remark}
\subsection{The Monte Carlo Euler method}\label{Kapote}
%
We now present
the Monte Carlo Euler time stepping
algorithm which will be the building block
for adaptive algorithms, with either quasi-deterministic
or stochastic time steps. For each realization
we first construct the jumps and its marks, and
conditioned on this information
we construct the approximate solution
$\barX$
as follows.
\par\noindent\smallskip\par\noindent
%
\fbox{
\begin{minipage}[h]{14cm}
%
%
{\tt Monte Carlo Euler time stepping algorithm}

\par\medskip
\textbf{Input:} Give a number $N+1$ of time nodes
$0={\tilde t}_0<\ttt_1<\dots<\ttt_N=T$,
and sample jump times
$0<\tau_1(\omega)<\dots<\tau_{\ssy \widehat{N}(\omega)}(\omega)< T$
with corresponding marks
$Z_1(\omega),\dots,Z_{\ssy \widehat{N}(\omega)}(\omega)$,
as explained in Remark \ref{Remarkk_1}.
\par\medskip
Set the jump counter $k=1$.
\par\medskip
\textbf{Time stepping:}
Consider an augmented partition given by the union
$$\{t_n(\omega)\}_{n=0}^{\ssy N_{\ssy A}(\omega)}
=\{\ttt_n\}_{n=0}^{\ssy N}\cup\{\tau_n(\omega)\}_{n=0}^{\ssy \widehat{N}(\omega)}$$
\hspace{0.5cm} where $N_{\ssy A}(\omega) =  N+\widehat{N}(\omega)$ ($a.s.$) 
is the number of time steps.

\hspace{0.5cm} Sample $X_0(\omega)$ and set the initial condition $\barX(t_{0},\omega)=X_0(\omega)$.

\hspace{0.5cm}\textbf{For} time steps $n = 0,\ldots,N_{\ssy A}(\omega)-1$

\hspace{1.0cm} Compute the remainder approximate grid values, $\barX(t_{n+1},\omega)$,

\hspace{1.0cm} by first constructing the left limit value of the approximated process,

\begin{equation}\label{eq:locdiffX}
\aligned
\barX(t_{n+1}^{-},\omega)&=\barX(t_n,\omega)+a(t_n,\barX(t_n,\omega))(t_{n+1}-t_n)\\
&\ken+\sum_{\ell = 1}^{\ell_0}b^{\ell}(t_n,\barX(t_n,\omega))(W^{\ell}(t_{n+1},\omega)-W^{\ell}(t_n,\omega)).
\endaligned
\end{equation}
\hspace{1.0cm}  When needed, introduce the correction due to jump discontinuities:

\hspace{1.0cm}\textbf{if} ($t_{n+1}=\tau_k(\omega)$) \textbf{then}

\hspace{1.5cm}
\begin{align}
\label{eq:jump}
&\barX(t_{n+1},\omega)=\barX(t_{n+1}^{-},\omega)+c\big(t_{n+1},\barX(t_{n+1}^{-},\omega),Z_{k}(\omega)\big),\\
&\text{increase $k$ to $k+1$}\notag
\end{align}

\hspace{1.5cm}\textbf{else}
\hspace{1.5cm}

\begin{equation}\label{eq:nojump}
\barX(t_{n+1},\omega)=\barX(t_{n+1}^{-},\omega),
\end{equation}

\hspace{1.0cm}\textbf{end-if}

\hspace{0.5cm}\textbf{end-for}

%
%
\end{minipage}
}

\par\noindent
\smallskip
\par\noindent
\par
When the initial time nodes
$\{\ttt_0,\dots,\ttt_{\ssy N}\}$
are the same for all realizations, we refer to
\emph{quasi-deterministic} time steps, or sometimes, simply \emph{deterministic}
time steps.
Otherwise, we speak about \emph{stochastic} time steps.

\par Besides, in the particular case when $a=b=0$
the approximate process $\barX$ has the same law
as $X$, due to the form of the grid proposed in
the numerical method (see \cite{glasserman}).
A similar situation occurs when only $b=0$: there,
conditioned to the realizations of ${\widehat N}$,
a higher order method for ODE integration
should be used to approximate $\barX$ between jumps.
It must be also noticed that the jump intensity $\lambda$
does not depend on the current value $X(t)$ of the process.
Such a dependence carries the necessity of implementing a different Monte Carlo Euler
algorithm, where the jumps and its marks are
sampled simultaneously with the trajectory of
the process, generating an aditional error
in the approximation scheme, as it does not seems
possible to have an exact simulation of the jump
structure in law, something that is possible in the present case.
%
\subsection{Error Control and Adaptivity}
%
The aim, for a given $\text{\rm TOL}>0$,
is to choose the size of time steps,
either quasi-deterministic or stochastic,
$$
\Delta{\ttt_n}= \ttt_{n+1}-\ttt_n,\quad n=0,\dots,N-1,
$$
and the number $M$ of
independent identically distributed samples
$\big\{\barX(\cdot\,,\omega_j)\big\}_{j=1}^{\ssy M}$
such that the computational work, 
defined as $M$ times the average of timesteps, i.e.
$M\times E[N_{\ssy A}]=M\times\big(E[N]+\Lambda(T)\big)$,
is minimal, constrained by the condition that the computational error
$\errorc$
defined by
\begin{equation}\label{eq:ce}
\errorc=E[g(X(T))]-\tfrac{1}{M}\sum_{j=1}^Mg(\barX(T,\omega_j)),
\end{equation}
is such that the event
\begin{equation}\label{eq:1.3}
\big|
\errorc
\big|\leq\tol
\end{equation}
has a  probability close to one.
The computational error
$\errorc$
naturally separates
as the sum of the (deterministic) time discretization error $\errort$ and the
statistical error $\errors$ given by
\begin{align}
\errort&=E\bigl[g(X(T))\bigr]-E\bigl[g(\barX(T))\bigr],\label{eq:errort}\\
\errors&=E[g(\barX(T))]-\tfrac{1}{M}\sum_{j=1}^Mg(\barX(T,\omega_j)).\label{eq:errors}
\end{align}
The time steps to construct the trajectories $\barX$ are determined from
statistical approximations of the time discretization error ${\me}_T$;
while the number $M$ of independent realizations of $\barX$,
is determined from  the statistical error ${\me}_S$.
Therefore, the number $M$ of realizations can be asymptotically determined
by standard limit theorems for sums of independent random variables.
\par
Efficient adaptive time stepping methods,
with theoretical basis,
use a posteriori error information, since
the a priori knowledge can usually not be
as precise as the a posteriori.
The present work develops adaptive time stepping methods
by proving  in Theorem  \ref{thm:2.2} and Theorem
\ref{thm:stoch_err_exp}
error estimates of ${\me}_T$
with leading order terms in computable a posteriori form.

The main
reference of Theorems \ref{thm:2.2} and
\ref{thm:stoch_err_exp}
is the work by Szepessy, Tempone and Zouraris \cite{STZ}
where similar techniques were applied considering the solution of
a stochastic differential equations, i.e. our case when $c= 0$.
The general inspiration of the presented results is the work by
Talay and Tubaro in \cite{TT} and its subsequent extension by
Protter and Talay  in \cite{PT} to stochastic differential equations
driven by L\'evy processes.

The main new idea here is the extension of the efficient use of
stochastic flows and dual functions to obtain the error expansion with
computable leading order term in Theorems \ref{thm:2.2} and
\ref{thm:stoch_err_exp}
in presence of jumps.
The use of dual functions is standard in optimal control theory
and in particular for adaptive mesh control
for ordinary and partial differential equations, see
\cite{BMV}, \cite{Jh}, \cite{JSz}, \cite{EEHJ}, and \cite{BR},
and was successfully applied in the probabilistic context in
\cite{STZ}, \cite{MSTZ}, \cite{DMSST} and in \cite{MSST}.

Concerning time steps, this work proposes
adaptive algorithms that  use  either
\begin{itemize}
\item ``quasi-deterministic''
time steps (see section \ref{sec:determ}),
in the sense that, although the grid is random, its
randomness  only depends
on the jumps, being in this sense is a strict generalization of
the deterministic time step algorithm in
\cite{STZ}.
\item ``stochastic'' time steps, meaning that the grid is random
also with respect to the Wiener measure, similarly as with the
stochastic time steps introduced in \cite{STZ}.
\end{itemize}
Theorem \ref{thm:2.2} describes the basis of an adaptive
algorithm  to estimate
the computational error;
the deterministic time steps are then chosen by
$$
(\Delta \ttt_n)^2 |E[\rho( \ttt_n,\omega)]|=\hbox{constant},
$$
-- observe that we can only control the deterministic points $ \ttt_n$
and not the whole augmented grid  $ \{t_n(\omega)\}_{n=0}^{\ssy N_{\ssy A}(\omega)}$ --
where $\rho(\ttt_n,\omega)$  is the function defined
in \eqref{eq:rhodef},
based on the weight functions $\varphi$ and $\varphi'$ defined
in \eqref{eq:phiend} -- \eqref{eq:phiinner}.
Provided the path $\barX(t_n), \ n=0,\ldots N_{\ssy A}$, is stored,
the leading order error bound can be evaluated
by solving, step by step, the two backward problems
\eqref{eq:phiend} -- \eqref{eq:phiinner}.
The backward evolutions
\eqref{eq:phiend} -- \eqref{eq:phiinner}
of the weight functions $\varphi$ and $\varphi'$
avoid solving for the
two variables $t,s$ present in $\frac{\partial\barX(t)}{\partial x(s)}
$, which
appears in the forward $t-$evolution equation for
$\frac{\partial\barX(t)}{\partial x(s)}$
in the identity $\varphi_i(t_n)=\partial_j g(\barX(T))
\partial\barX_j(t_n)/\partial x_i(t_n)$, cf. \eqref{eq:2.37}.  A
solution with two variables $s$ and $ t$ would
require work of the order $N^2$ for each realization, instead
of the corresponding work of the  order $N$ in Theorem~\ref
{thm:2.2}.

The second algorithm, presented in Section \ref{sec:sto},
is based on the expansion derived
in Theorem \ref{thm:stoch_err_exp} and
uses time steps
which may vary for different realizations
of the Wiener process.
Stochastic time steps are
advantageous for problems with singularities at random times.
The idea in this case is to choose the steps by
$$
(\Delta \ttt_n)^2 |\rho( \ttt_n,\omega)|=\hbox{constant},
$$
(where, in comparison with the previous algorithm,
there is no expectation), and this is achieved
through a test performed at each interval of each realization,
to decide whether to refine or not the given interval.
In this case,
when a node is added, the interpolation
is carried through the consideration of a Brownian
bridge to obtain the value of the approximated process
$\barX$ in the new added point.

Besides, since their use entails more work per realization than the
deterministic
time steps they  should be judiciously used. A natural application of
stochastic
time steps appears in the weak approximation of killed diffusions,
see \cite{DMSST}.
The optimal stochastic steps depend on the whole solution
$\barX(t), \ 0<t<T$, and in particular the step $\Delta t(t)$ at time
$t$
depends also on $W(u)$ for $t<u$. In stochastic analysis the concept of
\emph{adaptedness} means, intituively, that the values of the process at time $t$ depend
only on the events generated by the sources of randomness up to time $t$, i.e. do not
depend on future events.
In numerical analysis a method is said to be \emph{adaptive} when the approximate
solution is used to control the error, e.g. to determine the time steps.
Our stochastic time stepping algorithm is in this sense adaptive and non adapted, 
since the time steps $\Delta t(t)$ depends on values of $W(u)$ for $t<u$, i.e. on future
values of the Wiener process.
The stochastic time stepping algorithm (that achives higher precision)
requires, on one side, an aditional theoretical developement
including the introduction of Malliavin derivatives, and, on the other, 
requires the approximate computation of derivatives of second order. 
This second requirement is performed introducing the second derivatives $\varphi''$ of the fluxes
(the procedure is described in the Appendix). 
This computation is performed also in linear time.
%
%
%
%
%
%
\par
The  focus in the paper is on error estimates
for weak convergence of stochastic differential equations for diffusions with jumps.
Two particular important cases must be distinguished.
In first place, taking $c=0$ we obtain previous estimates in \cite{STZ}. 
More relevant in this instance is the situation when $a = b= 0$ and $c\neq 0$.
In this case, we are considering a pure jump process,
with finite number of jumps in the time interval $[0,T]$.
From the stochastic point of view, the situation is simpler,
and the process $\barX$ has exactly the same distribution as $X$.
This means that the time discretization error is null, and the problem
reduces to controlling the statistical error.

Furthermore, the deterministic problems associated with the computation
of the expected value $E[g(X(T))]$ via the  Feynman-Kac formula
are non local, in the sense that they involve integro-differential equations. 
Direct discretization of such time dependent integro-differential equations 
needs to approximate integrals at each time step, while the Euler Monte Carlo 
method avoids these expensive computations.

When $d$ is small and the Fourier transform of $X(t)$ is known
(for instance, $X$ is a L\'evy process), it is efficient to approximate $E[g(X(T))]$
based on Parseval's identity for this transform and the Fourier transform of $g$
\cite{Lewis,CT}.
The work \cite{AA} uses operator splitting to approximate the integral
term explicitly in Fourier space,
while approximating the other terms in the equation implicitly.
The work \cite{MPS} discretizes the partial integro-differential equation
by the $\theta$-scheme in time and a wavelet Galerkin method in space.
The resulting full Galerkin matrix is then replaced with a sparse matrix
in the wavelet basis, and the linear systems for each time step are solved
approximatively with GMRES in linear complexity.
The deterministic algorithm gives optimal convergence rates,
up to logarithmic terms, for the computed solution in the same complexity
as finite difference approximations of the standard Black-Scholes equation.
Other works include \cite{LiuLi}, where weak convergence schemes are analyzed,
\cite{higham} where two implicit numerical methods
for diffusion with jumps are presented,
analyzing strong convergence and nonlinear stability;
and \cite{bruti} that contains a survey, including some new results
on weak approximation for jump diffusion equations.
The technique used here is based on the 
Kolmogorov's backward equation developed in \cite{SV}
and \cite{SV1} to analyze uniqueness and
dependence on initial conditions for weak solutions of
stochastic differential equations with jumps.
The rest of the paper is as follows.
Section~\ref{sec:det_tstep} proves  error estimates for quasi-deterministic
time steps.
Section~\ref{Sec:StochTimStep} proves error estimates for stochastic time steps.
Section~\ref{sec:algo} presents implementations of adaptive algorithms
and, finally, Section~\ref{sec:numerical_experiments} includes results from
numerical experiments.
%
%
%
%
%
\section{An Error Estimate of the computational error
with deterministic time steps}\label{sec:det_tstep}
%
%
In this section we present in Theorem \ref{thm:2.2} 
an error expansion in a posteriori computable 
form for the time discretization error $\errort$.
The starting point for the analysis is
Lemma~\ref{lem:2.1} below. It uses the fact that the Euler method can be extended,
for theoretical purposes only, by
\begin{equation}\label{eq:2.1}
\begin{split}
\barX(t)-\barX(t_n)=&\int_{t_n}^t\brr{a}(s;\barX)ds
+\sum_{\ell=1}^{\ell_0}
\int_{t_n}^t\brr{b}^{\ell}(s;\barX)\,dW^{\ell}(s)
\quad\forall\,t\in\big[\,t_n,t_{n+1}\,\big),\\
\end{split}
\end{equation}
where $\brr{a}$ and $\brr{b}^\ell$ are the piecewise constant stochastic
approximations
\begin{equation}\label{eq:2.2}
\brr{a}(s;\barX)=a(t_n,\barX(t_n))
\ken\text{\rm and}\ken\brr{b}^\ell(s;\barX)=b^\ell(t_n,\barX(t_n))
\ken\text{\rm for}\ken s\in[t_n,t_{n+1}).
\end{equation}
Observe that the presence of jumps is realized in $\barX$ in
\eqref{eq:jump}, making possible not to introduce
a modified coefficient for $c$.
For simplicity we introduce the notation
\begin{eqnarray*}
\pa_k  \equiv \f{\pa}{\pa x_k},
\quad
\pa_{ki} \equiv \f{\pa^2}{\pa x_k \pa x_i},
\,
\dots
\,
\pa_t \equiv \f{\pa}{\pa t},
\end{eqnarray*}
and use the summation convention, i.e.,
if the same subscript appears twice in a term, the
term denotes the sum over the range of this subscript, e.g.
$c_{ik}\partial_k b_{j}\equiv\sum_{k=1}^d c_{ik}\partial_kb_{j}$,
and consequently
\begin{equation*}
d_{ij} \equiv 
\f{1}{2}\sum_{\ell=1}^{\ell_0} b_i^{\ell} b_j^{\ell}
=
\f{1}{2}b_i^{\ell} b_j^{\ell},
\quad
\bard_{ij}\equiv
\f{1}{2}\sum_{\ell=1}^{\ell_0}\barb_i^{\ell}\barb_j^{\ell}=
\f{1}{2}\barb_i^{\ell}\barb_j^{\ell}.
\end{equation*}
For a derivative $\partial_\alpha$ the notation $|\alpha|$ is
its order.
%
%
\begin{lemma}\label{lem:2.1}
Suppose that for some $m_0>[d/2]+10$
there are positive constants $k$ and $C$ such that
\begin{itemize}
\item[\rm (i)] $g\in{\mc}_{loc}^{m_0}(\rset^d)$ with
$|\partial_\alpha g(x)|\le C(1+|x|^k)$
for all $|\alpha|\le m_0$,
\item[\rm (ii)]
$E\big[\,|X(0)|^{2k+d+1}+|\barX(0)|^{2k+d+1}\,\big]\leq C$, and
\item[\rm (iii)]
$a$ and $b$ are bounded in ${\mc}^{m_0}([0,T]\times\rset^d)$
and the same holds uniformly for $c(\cdot,\cdot,z)$ for all $z\in\Z$.
\item[\rm (iv)] the initial data
$X(0)$ and its approximation $\barX(0)$ have the
same distribution.   
\end{itemize}
%
%
%
%
Then, the time discretization error in \eqref{eq:errort} satisfy 
\begin{equation}\label{eq:2.3}
\errort=\int_0^T E\big[\big( a_k(t,\barX(t^-))-\bara_k(t;\barX) \bigr)\partial_{k}u(t,\barX(t^-))
%
+\bigl(d_{ij}(t,\barX(t^-))-\bard_{ij}(t;\barX)\bigr)\partial_{ij}u(t,\barX(t^-))\bigr]dt
\end{equation}
where
\begin{equation}\label{eq:2.4}
u(t,x)=E\bigl[g(X(T))\mid X(t)=x\bigr].
\end{equation}
is the \emph{cost to go} function.
\end{lemma}
%
%
%
\begin{remark} We can relax condition {\rm (iv)} in the assumptions of the Lemma, and get an additional term 
of the form $E[u(0,X(0))-u(0,\barX(0))]$ in the error expansion in \eqref{eq:2.3}.
\end{remark}
%
%
%
%
%
%
\begin{proof}
There exists a unique solution $u\in {\mc}_{loc}^{1,6}([0,T]\times\rset^d)$
of the Kolmogorov backward equation
\begin{equation}\label{eq:2.12}
\begin{split}
\mathcal{L}^Xu(t,x)&\equiv{\pa_t}u(t,x)
+a_k\partial_{k}u(t,x)+d_{kn}\partial_{kn}u(t,x)
+\lambda(t)\,\int_{\Z}
[u(t,x+c(t,x,z))-u(t,x)]\mF(t,dz)=0\\
u(T,\cdot)&=g,\\
\end{split}
\end{equation}
satisfying the polynomial growth condition
\begin{equation*}
\max_{0\leq t\leq T}
|\partial_{\alpha}u(t,x)|\leq\,C\,\Big(\,1+|x|^{k+\frac{d+1}{2}}
\,\Big)\quad\forall\,|\alpha|\le 6,
\end{equation*}
for some positive $k$ and $C$ (cf. \cite{Fr}).
The Feynman-Kac formula without potential, 
implies that the solution $u$ 
of \eqref{eq:2.12} can be represented by the expected value in \eqref{eq:2.4}.
%
%
%
The $\Ito$ formula applied to \eqref{eq:2.1} (cf. \cite{IW}, p. 66) gives
$$\aligned
u(T,\barX(T))-u(0,\barX(0))
&=\int_0^T
\bigl\{\pa_tu(t,\barX(t^-))
+\bara_i(t;\barX )\partial_iu(t,\barX(t^-))
+\bard_{ij}\partial_{ij}u(t,\barX(t^-))\\
&\quad
+\lambda(t)\int_{\Z}[u(t,\barX(t^{-})+c(t,\barX(t^{-}),z))-u(t,\barX(t^{-}))]\mF(t,dz)
\bigr\}dt\\
&+\int_0^T \barb^\ell_i(t;\barX )\partial_i u(t,\barX(t^-)) dW^\ell(t)\\
&+\int_0^T\!\!\int_{\Z}[u(t,\barX(t^{-})+c(t,\barX(t^{-}),z))-u(t,\barX(t^{-}))]
\big(p(dt,dz)-q(dt,dz)\big)\\
\endaligned
$$
which, combined with \eqref{eq:2.12} to substitute 
$\pa_tu(t,\barX(t^-))$, yields
\begin{equation}\label{eq:2.14}
\begin{split}
u(0,\barX(0))-g(\barX(T))=&
\int_0^T\bigl(a_i(t,\barX(t^-))
- \bara_i(t;\barX )\bigr)\partial_i u(t,\barX(t^-))dt\\
&+\int_0^T\bigl(d_{ij}(t,\barX(t^-))-\bard_{ij}(t;\barX )\bigr)\partial_{ij}u(t,\barX(t^-))dt\\
&-\int_0^T \barb_i^\ell(t;\barX )\partial_i u(t,\barX(t^-))dW^\ell(t)\\
&\hskip-2.5truecm
-\int_0^T\!\!\int_{\ssy\Z}
\big[\,u(t,\barX(t^{-})+c(t,\barX(t^{-}),z))
-u(t,\barX(t^{-}))]\big(p(dt,dz)-q(dt,dz)\big).
\end{split}
\end{equation}
The expected value of the last two integrals is zero,
the first one by the martingale property of \Ito integrals,
and the second one due to the fact that
$p(dt,dz)-q(dt,dz)$
is a compensated random measure (i.e. a martingale measure).
Condition (iv) in the Lemma, and the representation of $u$ in \eqref{eq:2.4} show that
$$E
\big[u(0,\barX(0))\big]=E\big[u(0,X(0))\big] =E\big[g(X(T))\big].
$$
Therefore, taking expected values in both sides of
\eqref{eq:2.14} we arrive at the error representation
\eqref{eq:2.3}.
\end{proof}
%
%
%
%
\par
Lemma~\ref{lem:2.1} is combined with stochastic flows to  derive the a posteriori error expansion
in Theorem \ref{thm:2.2} below. 
This error expansion is based on the variations  of the processes 
$X$ and $\barX$.
For a process $\barX$, the first variation of a function $F(\barX(T))$
with respect to a perturbation in the initial location of the path $\barX$,
at time $s$, is denoted by
\begin{equation}\label{eq:2.15}
F'(T;s)=\partial_{x(s)}F(\barX(T))\equiv
\big(\partial_1F(\barX(T);\ \barX(s)=x),
\dots,
\partial_dF(\barX(T);\ \barX(s)=x)\big).
\end{equation}
The proof of Theorem \ref{thm:2.2} uses mainly that
the error in replacing  $g(X(T))$ in Lemma \ref{lem:2.1} by $g(\barX(T))$,
in the representation \eqref{eq:2.4} of $\partial_\alpha u$, yields
the small deterministic
remainder term $\int_0^TO\bigl((\Delta t)^2\bigr)dt$ in \eqref{eq:2.5}
of Theorem \ref{thm:2.2},
which is analogous to the $\OM(N^{-2})$
term in Talay and Tubaro's expansion, cf. \cite{TT},
 and needs some a priori estimate
to be controlled. Lemma \ref{lem:2.1} can be  applied
to estimate this error.
The second important ingredient in the proof is
the Markov property of $\barX$ satisfied at the discrete times $t_n$.
Based on the fact that $\barX(t_n)$ is ${\cF}_{t_n}$ measurable,
the nested expected values
$$
E\bigl[ a_j(t_n,\barX(t_n))
\partial_{x_j(t_n)} E[g(\barX(T))\mid{\cF}_{t_n}] \bigr]
$$
in \eqref{eq:2.3} can,
by  the definition
 of $\varphi$ and its implication \eqref{eq:2.37}, be decoupled to
$$
E\bigl[a_j(t_n,\barX(t_n)) \varphi_j(t_n)\bigl],
$$
which reduces the computational complexity substantially,
see Lemma~\ref{lem:2.5}.
\begin{thm}[Error expansion with deterministic time steps]\label{thm:2.2}
Suppose that $a$, $b$, $g$, $X$ and $\barX$,  satisfy the
assumptions
in Lemma \ref{lem:2.1}. 
Then, the time discretization error in \eqref{eq:errort} has the expansion
\begin{equation}\label{eq:2.5}
\errort
=
E\left[\sum_{\ttn=0}^{N-1}\rho(\ttt_\ttn,\omega)(\Delta \ttt_\ttn)^2\right]
+E\left[\sum_{\ttn=0}^{N-1} (\Delta \ttt_\ttn)^2
\Bigl\{\OM(\Delta \ttt_\ttn)+
        \sum_{m=n}^{N-1} \OM((\Delta \ttt_{\ttn})^2)\Bigr\}\right]
\end{equation}
where the leading order error term is in computable a posteriori form.
\begin{multline}
\label{eq:rhodef}
\rho(\ttt_\ttn,\omega) \equiv
\half
\sum_{n \in \cJ_\ttn}\bigl[
                 \bigl(a_i(t_{n+1},\barX(t_{n+1}^{-},\omega)) -
                       a_i(t_{n},\barX(t_{n},\omega))\bigr)\varphi_i(t_{n+1}^{-},\omega)
                     \bigr]
\frac{\Delta t_n}{(\Delta \ttt_\ttn)^2}\\
+\half
\sum_{n \in \cJ_\ttn  }\bigl[
                 \bigl(d_{ik}(t_{n+1},\barX(t_{n+1}^{-},\omega)) -
                       d_{ik}(t_{n},\barX(t_{n},\omega))\bigr)\varphi'_{ik}(t_{n+1}^{-},\omega)
                     \bigr]
\frac{\Delta t_n}{(\Delta \ttt_\ttn)^2}
\end{multline}
with $\cJ_\ttn \equiv \{n: \ttt_\ttn \le t_n<\ttt_{\ttn+1} \}$, $\ttn = 0,\ldots,N-1$,
and based on the discrete dual functions $\varphi(t_n)\in \rset^d$
and $\varphi'(t_n)\in \rset^{d\times d}$, which are determined as follows.
The function $\varphi$ and and its first variation
\begin{equation}\label{eq:2.8}
\varphi'_{ik}(t_n,\omega)=\partial_{x_k(t_n)} \varphi_i(t_n,\omega)
\equiv {\frac{\partial\varphi_i(t_n; \barX(t_n)=x)}{\partial x_k}}
\end{equation}
satisfies \eqref{eq:phiend} -- \eqref{eq:phiinner}.
\end{thm}
\begin{remark}
When implementing an algorithm based on this result the expectation of the
sum of errors in \eqref{eq:2.5} is approximated by the mean of the 
errors along the $M$ simulated trajectories, i.e.:
$$
E\left[\sum_{\ttn=0}^{N-1}\rho(\ttt_\ttn,\omega)(\Delta \ttt_\ttn)^2\right]
\sim
\frac 1M\sum_{j=1}^M\sum_{\ttn=0}^{N-1}\frac{\rho(\ttt_\ttn,\omega_j)}{M}{(\Delta \ttt_\ttn)^2}.
$$
The statistical error of this approximation can be expressed as
$$
E\left[\sum_{\ttn=0}^{N-1}\rho(\ttt_\ttn,\omega)(\Delta \ttt_\ttn)^2\right]
-
\frac 1M\sum_{j=1}^M\sum_{\ttn=0}^{N-1}\frac{\rho(\ttt_\ttn,\omega_j)}{M}{(\Delta \ttt_\ttn)^2}
=\int_0^T( I_M + II_M) dt
$$
where the distributions of the statistical errors $\sqrt M I_M$ and $\sqrt M II_M$
weakly converge to normal distributions with mean zero and time interval dependent 
variances given by
$$
\var \bigl[
\sum_{n \in \cJ_\ttn}
(a_i(t_{n+1},\barX(t_{n+1}^{-})) -
a_i(t_{n},\barX(t_{n}))\, )\varphi_i(t_{n+1}^{-})
\bigr]=\OM(\Delta \ttt_\ttn),
$$
and
$$
\var \bigl[
\sum_{n \in \cJ_\ttn}
(d_{ik}(t_{n+1},\barX(t_{n+1}^{-})) -
d_{ik}(t_{n},\barX(t_{n}))\, )\varphi'_{ik}(t_{n+1}^{-})\bigr]
=\OM(\Delta \ttt_\ttn),\\
$$
respectively.
\end{remark}
\begin{proof}[Proof of Theorem \ref{thm:2.2}]
The main content of Theorem \ref{thm:2.2} is the replacement of the (non computable) 
estimate in Lemma \ref{lem:2.1}
by an expansion with computable leading order term. 
For this purpose the derivatives of the expected value 
$$
\partial_\alpha u(x,t)=\partial_\alpha E[g(X(T))\mid X(t)=x]
$$
appearing in the integral in Lemma \ref{lem:2.1} 
are approximated by  the corresponding derivatives of
$$
\bar u(x,t)\equiv E[g(\barX(T))\mid \barX(t)=x],
$$
that depends on the simulated solution $\barX$.
The proof 
is divided into three steps:
\begin{itemize}
\item[(i)] in Lemma \ref{lem:2.3} we estimate the quadrature error;
\item[(ii)] in Lemma \ref{lem:2.4} we bound the error in replacing
$\partial_\alpha u$ by $\partial_\alpha \bar u$ with the use of stochastic flows and 
its variations;
\item[(iii)] in Lemma \ref{lem:2.5} we use the discrete dual functions 
$\varphi$ and $\varphi'$ (that solve the backward evolution problems
see \eqref{eq:phiend} -- \eqref{eq:phiinner} in the Appendix) to 
derive a computable representation of $\partial_\alpha \bar u$.
\end{itemize}

We begin with the first step.
\begin{lemma}[Quadrature approximation]\label{lem:2.3}
Suppose that the assumptions in Lemma~\ref{lem:2.1} hold.
Let $\cJ_\ttn \equiv \{n: \ttt_\ttn \le t_n<\ttt_{\ttn+1} \}$,
$\ttn = 0,\dots,N-1$.
Then the quadrature error terms satisfy
$$
\aligned
&\int_{\ttt_\ttn}^{\ttt_{\ttn+1}}E\bigl[\bigl(a_{i}(t,\barX(t^-))-\bara_{i}(t;\barX)
\bigr)
\partial_{i}u(t,\barX(t^-))\bigr] dt\\
&- E\bigl[\sum_{n \in \cJ_\ttn}
         \bigl(a_{i}(t_{n+1},\barX(t_{n+1}^{-}))-a_{i}(t_n,\barX(t_n) )
\bigr)
\partial_{i}u(t_{n+1},\barX(t_{n+1}^{-}))\f{\Delta t_n}{2} \bigr]
= O\bigl((\Delta \ttt_\ttn)^3\bigr),\\
\endaligned
$$
and
$$
\aligned
&\int_{\ttt_\ttn}^{\ttt_{\ttn+1}}
E\bigl[\bigl(d_{ij}(t,\barX(t^-))-\bard_{ij}(t;\barX)
\bigr)
\partial_{ij}u(t,\barX(t^-))\bigr] dt\\
&- E\bigl[\sum_{n \in \cJ_\ttn}
          \bigl(d_{ij}(t_{n+1},\barX(t_{n+1}^{-}))-d_{ij}(t_n,\barX(t_n) )
\bigr)
\partial_{ij}u(t_{n+1},\barX(t_{n+1}^{-}))\f{\Delta t_n}{2}\bigr]
= O\bigl((\Delta \ttt_\ttn)^3\bigr).\\
\endaligned
$$
\end{lemma}
\begin{proof}
Denote by $\sigmatt$ the $\sigma$-algebra generated by the jumps and marks
in $[0,T]$ constructed in Remark \ref{Remarkk_1}.
Then 
%
\begin{multline}\label{eq:condexquad}
\int_{\ttt_\ttn}^{\ttt_{\ttn+1}}
E\bigl[\bigl(d_{ij}(t,\barX(t^-))-\bard_{ij}(t;\barX)
\bigr)
\partial_{ij}u(t,\barX(t^-))\bigr] dt\\
 = E\bigl[\sum_{n \in \cJ_\ttn}
 \int_{t_n}^{t_{n+1}}
  E\bigl[
  \bigl(
      d_{ij}(t,\barX(t^-))-\bard_{ij}(t;\barX)
  \bigr)
  \partial_{ij}u(t,\barX(t^-))
  \mid\sigmatt \bigr]
 dt \bigr].
\end{multline}
Observe that in $[t_n,t_{n+1})$ the conditioned
process $\barX$ has no jump discontinuities, 
and introduce the notations
$$
\aligned
h(t,\barX(t^-))&\equiv \bigl(d_{ij}(t,\barX(t^-))-\bard_{ij}(t;\barX)\bigr)
\partial_{ij}u(t,\barX(t^-)),\\
\bar h(t)&\equiv \f{t-t_n}{\Delta t_n}
\bigl(d_{ij}(t_{n+1},\barX(t_{n+1}^{-}))-\bard_{ij}(t_n;\barX )\bigr)
\partial_{ij}u(t_{n+1},\barX(t_{n+1}^{-})).\\
\endaligned
$$
Then the quadrature error  satisfies
\begin{multline*}
\int_{t_n}^{t_{n+1}}E[h(t,\barX(t^-))-\bar h(t) \mid\sigmatt]dt
=\int_{t_n}^{t_{n+1}}{E}\bigl[\bigl(d_{ij}(t,\barX(t^-))-\bard_{ij}(t;\barX)\bigr)
\partial_{ij}u(t,\barX(t^-))\mid\sigmatt\bigr] dt\\
-E\bigl[\bigl(d_{ij}(t_{n+1},\barX(t_{n+1}^{-}))-\bard_{ij}(t_n;\barX )\bigr)
\partial_{ij}u(t_{n+1},\barX(t_{n+1}^{-})) \mid\sigmatt\bigr] \f{\Delta t_n}{2},
\end{multline*}
and 
$E[\bar h(t)\mid\sigmatt]$ 
is the linear nodal projection of the smooth function
$E[h(t,\barX(t^-))\mid\sigmatt]$ 
in the interval 
$[t_n,t_{n+1})$. 
Therefore, a standard interpolation estimate yields
\begin{equation}
\Big|\int_{t_n}^{t_{n+1}}E[h(t,\barX(t^-))-\bar h(t) \mid\sigmatt]dt\Big| 
\le
\f{1}{8}(\Delta t_n)^2
\int_{t_n}^{t_{n+1}}\bigl|\f{d^2}{dt^2}E[h(t,\barX(t^-))\mid\sigmatt]\bigr|dt. 
\label{eq:2.28}
\end{equation}
Denoting
$
\mathcal{L}h\equiv \pa_t h+\bara_i\partial_i h+\bard_{ij}\partial_{ij}h,
$
$\Ito$'s formula and condition (iii) in Lemma \ref{lem:2.1}
show that there exist constants $C_1,C_2$ such that, for $t \in (t_n,t_{n+1})$
\begin{equation}
\aligned
\label{eq:2.29}
\f{d}{dt}E[h(t,\barX(t^-)) \mid\sigmatt] &= E[\mathcal{L}h(t,\barX(t^-)) \mid\sigmatt]\le C_1,\\
\f{d^2}{dt^2}E[h(t,\barX(t^-)) \mid\sigmatt] &= E[\mathcal{L}^2h(t,\barX(t^-)) \mid\sigmatt]\le C_2,\\
\endaligned
\end{equation}
which combined with \eqref{eq:condexquad} and \eqref{eq:2.28}
proves the estimate of the diffusion term in the lemma. 
The estimate of the drift term follows analogously.
\end{proof}
In the second step of the proof the derivative $\partial_\alpha u$ and its approximation 
$\partial_\alpha\bar u$
are evaluated respectively through expected values of \emph{stochastic flows} of $X$ and $\barX$ and 
its \emph{variations}, that we now introduce.
%
Recall the definition \eqref{eq:2.15} of the first variation, and let
$$
\delta_{ik}
\equiv
\left\{
\aligned &0  \ \ i\ne k,\\
         &1 \  \ i= k.  \\
\endaligned
\right.
$$
The following equation for the first variation of the process $X$ at time $s>t$ hold:
\begin{equation}\label{eq:2.16}
\aligned
dX_{ij}'(s)&=\partial_{k}a_i(s,X(s^-)) X'_{kj}(s^-)ds
+\partial_{k}b^{\ell}_i(s,X(s^-))
X'_{kj}(s^-)dW^{\ell}(s)\\
&\hskip1.0truecm+\int_{\Z}\partial_{k}c_i(s,X(s^{-}),z)X'_{kj}(s^-)\,p(ds,dz),\\
X'_{ij}(t)&=\delta_{ij}.\\
\endaligned
\end{equation}
Similarly, for the second variation of the process $X$ at time $s>t$ we have:
\begin{equation}\label{eq:2.17}
\aligned
dX_{ijn}''(s)&=\bigl[\partial_{k}a_i(s,X(s^-)) X_{kjn}''(s^-)
+\partial_{kr}a_i(s,X(s^-))X'_{kj}(s^-)X'_{rn}(s^-)\bigr]ds\\
&\quad+\bigl[\partial_{k}b_i^{\ell} (s,X(s^-))X_{kjn}''(s^-)
+\partial_{kr}b_i^{\ell}(s,X(s^-))X_{kj}'(s^-)X_{rn}'(s^-)\bigr]dW^{\ell}(s)\\
&\quad+\int_{\Z}
\bigl[\partial_{k}c_i(s,X(s^-),z) X_{kjn}''(s^-)
+\partial_{kr}c_i(s,X(s^-),z)X'_{kj}(s^-)X'_{rn}(s^-)\bigr]p(ds,dz),\\
X''_{ijn}(t)&=0.\\
\endaligned
\end{equation}
For the third variation of the process $X$ at time $s>t$ we have:
\begin{equation}
\aligned
dX'''_{ijnm}(s)&= \bigl[\qqq_{k}a_i(s,X(s^{-})) X_{kjnm}'''(s^{-})+\qqq_{kr}a_i(s,X(s^-))X'_{kj}(s^-)X''_{rnm}(s^-)\\
        &\quad+\qqq_{kr}a_i(s,X(s^-))X'_{kn}(s^-)X''_{rjm}(s^-)+\qqq_{kr}a_i(s,X(s^-))X'_{km}(s^-)X''_{rjn}(s^-)\\
        &\quad+\qqq_{krv}a_i(s,X(s^-))X'_{kj}(s^-)X'_{rn}(s^-)X'_{vm}(s^-)\bigr]ds\\
&+\bigl[\qqq_{k}b^\ell_i(s,X(s^-)) X_{kjnm}'''(s^-)+\qqq_{kr}b^\ell_i(s,X(s^-))X'_{kj}(s^-)X''_{rnm}(s^-)\\
&\quad+\qqq_{kr}b^\ell_i(s,X(s^-))X'_{kn}(s^-)X''_{rjm}(s^-)+\qqq_{kr}b^\ell_i(s,X(s^-))X'_{km}(s^-)X''_{rjn}(s^-)\\
&\quad+\qqq_{krv}b^\ell_i(s,X(s^-))X'_{kj}(s^-)X'_{rn}(s^-)X'_{vm}(s^-)\bigr]dW^\ell(s)\\
&+\int_{\Z}\bigl[\qqq_{k}c_i(s,X(s^-)) X_{kjnm}'''(s^-)+\qqq_{kr}c_i(s,X(s^-),z)X'_{kj}(s^-)X''_{rnm}(s^-)\\
&\quad+\qqq_{kr}c_i(s,X(s^-),z)X'_{kn}(s^-)X''_{rjm}(s^-)+\qqq_{kr}c_i(s,X(s^-),z)X'_{km}(s^-)X''_{rjn}(s^-)\\
&\quad+\qqq_{krv}c_i(s,X(s^-),z)X'_{kj}(s^-)X'_{rn}(s^-)X'_{vm}(s^-)\bigr]p(ds,dz),\\
X'''_{ijnm}(t)&=0,
\endaligned\label{eq:2.18}
\end{equation}
%
and similarly for the fourth variation of the process $X$ at time $s>t$:
\begin{equation}\label{eq:2.19}
dX''''_{ijnmp}=\dots,
\qquad
X''''_{ijnmp}(t)=0.
\end{equation}
This equations imply the representation of the derivatives of expectations 
with  stochastic flows as follows (cf. \cite{Pro} and \cite{S}). 
For the first derivatives:
\begin{equation}
\qqq_{k}u(t,x)=E\bigl[\qqq_{i}g(X(T))
X'_{ik}(T)\mid  X'_{ij}(t)=\delta_{ij}, X(t)=x\bigr],\label{eq:2.20}
\end{equation}
for the second derivatives:
\begin{equation}
\aligned
\qqq_{kn}u(t,x)=E\bigl[\qqq_{i}g(X(T)) X''_{ikn}(T)
+
\qqq_{ir}g(X(T)) X_{ik}'(T) X_{rn}'(T)
\mid  X_{ikn}''(t)=0,
X'_{ij}(t)=\delta_{ij}, X(t)=x],\\
\endaligned\label{eq:2.21}
\end{equation}
for the third derivatives:
\begin{multline}
\qqq_{knm}u(t,x)=E\bigl[\qqq_{i}g(X(T)) X'''_{iknm}(T)+\qqq_{ir}g(X(T)) X_{ik}'(T) X_{rnm}''(T)\\
+\qqq_{ir}g(X(T)) X_{in}'(T) X_{rkm}''(T)+\qqq_{ir}g(X(T)) X_{im}'(T) X_{rkn}''(T)\\
+\qqq_{irv}g(X(T)) X_{ik}'(T) X_{rn}'(T)X'_{vm}(T)\mid  X_{iknm}'''(t)=X_{ikn}''(t)=0,X'_{ij}(t)=\delta_{ij}, X(t)=x\bigr],
\label{eq:2.22}
\end{multline}
and for the fourths derivatives:
\begin{equation}\label{eq:2.23}
\qqq_{knmp}u(t,x)=\dots \ \ .
\end{equation}
Let $Y=(X,X',X'', X''', X'''')^T$ and
let $I$ denote the $d\times d$ identity matrix. Then
the system (\ref{eq:2.16}-\ref{eq:2.19}) can be written
\begin{equation}
\aligned
dY       &=A(t,Y(t^-))dt
                 +B^{\ell}(t,Y(t^-))dW^{\ell}(t)+\int_{\Z}C(t,Y(t^-),z)p(dt,dz),\ken t>t_0\\
Y(t_0)   &=(x,I,0, 0, 0)^T.\\
\endaligned\label{eq:2.24}
\end{equation}
Furthermore, rewrite the representation (\ref{eq:2.20}-\ref{eq:2.23}) as
\begin{equation}
\aligned
f_i(Y)&\equiv\qqq_{k}g(X)X_{ki}',\\
f_{ij}(Y)&\equiv\qqq_{k}g(X)X_{kij}''
+\qqq_{kn}g(X)X_{ki}'X_{nj}',\\
f_{ijm}(Y)&\equiv\qqq_{k}g(X)X_{kijm}'''
+\qqq_{kn}g(X)X_{ki}'X_{njm}''\\
&\quad+\qqq_{kn}g(X)X_{kj}'X_{nim}''
+\qqq_{kn}g(X)X_{km}'X_{nji}''\\
&\quad +\qqq_{knv}g(X)X_{ki}'X_{nj}'X_{vm}',\\
f_{ijmn}(Y)&\equiv\ldots \ \ .\\
\endaligned\label{eq:2.25}
\end{equation}
The Euler approximation of $Y$, the solution \eqref{eq:2.24}, is denoted by 
$\barY=\big(\,\barY^0,\barY^1,\barY^2,\barY^3,\barY^4\,\big)^{\ssy T}$
and can be extended as the solution of the stochastic differential equation 
with piecewise constant drift and diffusion fluxes
\begin{equation}\label{eq:2.26}
d\barY=\brr{A}(t;\barY)dt+\brr{B}^{\ell}(t;\barY)dW^{\ell}(t)
+\int_{\ssy\Z}C(t,\barY,z)p(dt,dz),
\end{equation}
defined as in \eqref{eq:locdiffX}, \eqref{eq:nojump}
and \eqref{eq:2.1}, with coefficients defined as in \eqref{eq:2.2}.

An important consequence of the Euler method is that
the variation and the Euler discretization commute.
This yields for each $\alpha$ 
the representation
\begin{equation}\label{eq:alpha}
\partial_\alpha u(t,x)
-
\partial_\alpha\bar u(t,x)=
E\bigl[f_\alpha(Y(T))
-f_\alpha(\barY(T))\mid Y(t)
=\barY(t)=(x,I,0,0,0)^T\bigr].
\end{equation}
Now we rely on Lemma \ref{lem:2.1} applied to the process $Y$ to obtain the representation
\begin{equation}\label{eq:2.27}
\begin{split}
E\bigl[f_{\alpha}(Y(T))-&f_{\alpha}(\barY(T))\mid Y(t)=\barY(t)=(x,I,0,0,0)^T\bigr]\\
=&\int_t^TE\Big[
(A-\brr{A})_k\qqq_{k}\nu^{\alpha}(s,\barY(s^-))
+(D-\brr{D})_{kn}\qqq_{kn}\nu^{\alpha}(s,\barY(s^-))
\mid \barY(t)=(x,I,0,0,0)^T\Big]ds\\
\end{split}
\end{equation}
to be used in Lemmas \ref{lem:2.4} and \ref{lem:2.5} below,
where the cost to go function $\nu^{\alpha}$ is, for each $f_{\alpha}$
in \eqref{eq:2.25}, defined as
\begin{equation*}
\nu^{\alpha}(t,y)=E\bigl[f_{\alpha}(Y(T))\mid Y(t)=y\bigr],
\end{equation*}
we use the notation
$$
D_{kn}=\tfrac{1}{2} B_k^{\ell} B_n^{\ell},
\quad
\brr{D}_{kn}=\tfrac{1}{2}\brr{B}_k^{\ell}\brr{B}_n^{\ell},
$$
and the corresponding Kolmogorov Backward equation is
$$
\aligned
\mathcal{L}^Y\nu^{\alpha}&\equiv{\pa_t}\nu^{\alpha}+A_k\qqq_{k}\nu^{\alpha}
+D_{kn}\qqq_{kn}\nu^{\alpha}
+\lambda(t)\int_{\Z}[\nu^{\alpha}(t,Y+C(t,Y,z))-\nu^{\alpha}(t,Y)]\mF(t,dz)=0, \ \ t<T,\\
\nu^{\alpha}(T,\cdot)&=f_{\alpha}.\\
\endaligned
$$
We are ready for the second step.
%

%
\begin{lemma}[Approximation of $\partial_\alpha u$]\label{lem:2.4}
%
Let the piecewise constant mesh function $\Delta t$ be defined by
$$\Delta{t}(s)\equiv\Delta t_n
\ken\text{\rm for}\ken s\in [t_n,t_{n+1})
\ken\text{\rm and}\ken n=0,\dots,N_{\ssy A}(\omega)-1.$$
Suppose that the assumptions in Lemma \ref{lem:2.1} hold.
Then the discretization errors of the stochastic
flows, for $|\alpha|\le 4$
, satisfy
\begin{equation}
\partial_\alpha (u-\bar u)(t_n,\barX(t_n)) = \int_{t_n}^T E[\OM(\Delta t(s))\mid\cF_{t_n}]ds \\
=  
\OM(\Delta \ttt_{max}).
\label{eq:2.30}
\end{equation}
Furthermore
\begin{equation}\label{eq:2.31}
\aligned
E\bigl[\bigl(a_{i}(t_{n+1},\barX(t_{n+1}^{-}) ) -\bara_{i}(t_n;\barX
)\bigr)& \bigl(\partial_{i}u(t_{n+1},\barX(t_{n+1}^{-}))-
\partial_{i}\bar u(t_{n+1},\barX(t_{n+1}^{-}))\bigr)\mid\sigmatt\bigr] \\
&= \Delta t_n \int_{t_{n+1}}^T \OM(\Delta t(s)) ds,\\
\endaligned
\end{equation}
\begin{equation}
\label{eq:2.32}
\aligned
E\bigl[\bigl(d_{ij}(t_{n+1},\barX(t_{n+1}^{-}) )
-\bard_{ij}(t_n;\barX )\bigr)&
\bigl(\partial_{ij}u(t_{n+1},\barX(t_{n+1}^{-}))-
\partial_{ij}\bar u(t_{n+1},\barX(t_{n+1}^{-}))\bigr)\mid\sigmatt\bigr]\\
&= \Delta t_n \int_{t_{n+1}}^T \OM(\Delta t(s)) ds. \\
\endaligned
\end{equation}
\end{lemma}
\begin{proof}
The combination of \eqref{eq:alpha} and \eqref{eq:2.27} give
%
\begin{equation}\label{eq:2.33}
\begin{split}
\partial_\alpha(u-\bar u)(t_n,\barX(t_n))=\int_{t_n}^T E\bigl[(A_i-\brr{A}_i)\partial_i\nu^\alpha(s,\barY(s^-))
+(D_{ij}-\brr{D}_{ij})
\partial_{ij}\nu^\alpha(s,\barY(s^-))\mid {\cF}_{t_n}\bigr] ds
\end{split}
\end{equation}
Now, for $t_m\le t<t_{m+1}$, introduce the notation
$$
h(t,\barY(t))\equiv (A_i-\brr{A}_i)\partial_i\nu^\alpha(t,\barY(t))
+(D_{ij}-\brr{D}_{ij})\partial_{ij}\nu^\alpha(t,\barY(t)),
$$
and
\begin{align*}
\mathcal{L}_0w(t,y)        &\equiv \bigl(\pa_t w
+\brr{A}_n\partial_n w + \brr{D}_{kn}\partial_{kn}w \bigr)(t,y)\\
\mathcal{L}^{\barY}w(t,y)  &\equiv \mathcal{L}_0w(t,y)
+\lambda(t)\int_{\Z}[w(t,y+C(t,y,z))-w(t,y)]\mF(t,dz).
\end{align*}
Observe that $h(t_m,\barY(t_m))=0$, apply
$\Ito$'s formula in the interval $[t_m,t]$ to obtain
$$
E[h(t,\barY(t^-))\mid  {\cF}_{t_m}]=
\int_{t_m}^tE[\mathcal{L}^{\barY}h(s,\barY(s^-))\mid{\cF}_{t_m}]ds=\OM(\Delta t_m), \ \ t_m\le t < t_{m+1},
$$
where the bound follows as in \eqref{eq:2.29}. This plugged into
\eqref{eq:2.33} proves \eqref{eq:2.30}.

The estimate \eqref{eq:2.32} follows similarly by  defining now
$$
\tilde h(t,\barX(t))\equiv
(d_{ij}-\brr{d}_{ij})\partial_{ij}(u-\bar u)(t,\barX(t)).
$$
Then the $\Ito$'s formula shows as in (\ref{eq:2.28}-\ref{eq:2.29})
$$
E[\tilde h(t,\barX(t))\mid\sigmatt]=\int_{t_n}^{t}
E[\mathcal{L}_0\tilde{h}(s,\barX(s^-))\mid\sigmatt] ds,
\ \ t_n\le t < t_{n+1}.
$$
The final step to prove \eqref{eq:2.32} is to establish
$$
E[\mathcal{L}_0\,\tilde{h}(s,\barX(s^-))\mid\sigmatt]=\int_{s}^T\OM(\Delta t(\tau))d\tau.
$$
The function $\mathcal{L}_0\,\tilde h(s,\barX(s^-))$ splits into the two types of terms
$h_1\equiv (d-\brr{d})v$ and $h_2\equiv v\partial_\alpha(u-\bar u)$, with
smooth functions $v$ of $(s,\barX(s))$. The $\Ito$ formula  again
shows that
$$
E[h_1(s,\barX(s^-))\mid\sigmatt]= \int_{t_n}^s E[\mathcal{L}_0h_1(\tau,\barX(\tau^-))\mid\sigmatt] d\tau
=\OM(\Delta t_n), \ \ t_n\le s < t_{n+1}.
$$
Moreover \eqref{eq:2.30} implies
$$
E[h_2(s,\barX(s^-))\mid\sigmatt]= \int_{s}^T\OM(\Delta t(\tau))d\tau,
$$
and consequently \eqref{eq:2.32} holds.
The estimate \eqref{eq:2.31} of the drift terms follows analogously.
\end{proof}
%
\begin{defin}[Local discrete solution operator]\label{def:ldso}
Write the Euler time stepping for the time nodes
$t = t_0,t_1^-, t_{1},t_2^-, t_{2},\ldots$
as
\begin{equation}\label{eq:locsolop}
\barX(t_{next},\omega) = \nA(t,\barX(t),\omega)
\end{equation}
where the next computation time is given by
$$
t_{next}= 
\left\{
       \begin{array}{lll}
                     t_{n} & if & t = t_{n}^{-} \\
                     t_{n+1}^{-} & if & t = t_{n}.
       \end{array}
\right.
$$
Observe that depending on $\omega$ and $t$, 
we may have 
$$
\nA(t,x,\omega)=
\begin{cases}
x+a(t,x)\Delta{t}+b^\ell(t,x)\Delta{W^\ell},&\text{ if \ $t=t_n$},\\
x+ c(t, x,Z_k(\omega)),& \text{ if \ $t=t_{n}^{-}$ and $t_n$ is the $k$-th jump time},\\
x&\text{ if \ $t=t_{n}^{-}$ and $t_n$ is a not a jump time}.
\end{cases}
$$
\end{defin}

\begin{lemma}[Representation with discrete duals]
\label{lem:2.5}
Suppose that the assumptions in Lemma
\ref{lem:2.1} hold. Then the dual functions $\varphi$ and $\varphi'$,
defined Appendix \ref{app:dualeqs},
 satisfy for $t = t_n$ and  $t = t_{n+1}^{-}$
\begin{equation}\label{eq:2.34}
\aligned
\partial_i\bar u(t,\barX(t)) &= E[\varphi_i(t)\mid {\cF_{t}}],
\\
\endaligned
\end{equation}
\begin{equation}\label{eq:2.35}
\aligned
\partial_{ij}\bar u(t,\barX(t)) &= E[\varphi'_{ij}(t)\mid  {\cF_{t}}]. \\
\endaligned
\end{equation}
%
\end{lemma}
\begin{proof}
Equations \eqref{eq:2.1}, \eqref{eq:2.16} and \eqref{eq:2.20} show that the first
variation of the Euler approximation $\barX$ is in fact equal
to the Euler approximation of the first variation $X'$ and consequently
$$
\partial_i\bar u(t,\barX(t))=E[\partial_jg(\barX(T))\barX'_{ji}(T;t)\mid  {\cF_t}],
$$
where $\barX'_{ji}(s;t)\ (s>t)$, is the Euler approximation \eqref{eq:2.26}
of $X'$ with initial data $\barX'_{ji}(t;t)=\delta_{ji}$.

Let, for $t = \ldots,t_n, t_{n+1}^{-},\ldots$,
$$
\varphi_i(t)\equiv \partial_{x_i(t)} g(\barX(T)),
$$
i.e.
\begin{equation}\label{eq:2.37}
\varphi_i(t) = \partial_jg(\barX(T))\barX'_{ji}(T;t).
\end{equation}
We prove inductively that $\varphi_i(t)$ is
the solution of the corresponding problem in
\eqref{eq:phiend}-\eqref{eq:phiinner}.
Since \eqref{eq:phiend} is trivially true it remains to prove the
inductive step.
By the chain rule we have
$$
\aligned
\partial_kg(\barX(T))\barX'_{ki}(T;t) =&
\partial_kg(\barX(T))\barX'_{kj}(T;t_{next}) \barX'_{ji}(t_{next};t)
\endaligned
$$
or in other words
\begin{equation}\label{eq:phievol}
\varphi_i(t) = \varphi_j(t_{next}) \partial_i \nA_j(t,\barX(t))
\end{equation}
which is  equivalent to
\eqref{eq:phiend}-\eqref{eq:phiinner},
what we wanted to prove.

The equality \eqref{eq:2.37} implies that
$$
\partial_{ij}\bar u(t,\barX(t))
=E\bigl[\partial_{x_j(t)}\varphi_{i}(t)\mid  {\cF_{t}}].
$$
The next step is to verify that the first variation of $\varphi$,
\begin{equation}\label{eq:2.38}
\varphi'_{ij}(t)\equiv\partial_{x_j(t)}\varphi_i(t) \equiv
\frac{\partial\varphi_i(t; \barX(t)=x)}{\partial x_j},
\end{equation}
satisfies the backward recursive equation
\eqref{eq:phiend}-\eqref{eq:phiinner}.
First, differentiate the equation
\eqref{eq:phievol}
to obtain
\begin{equation}\label{eq:2.39}
\aligned \varphi'_{ik}(t)&=\partial_i {\nA}_j\partial_{x_k(t)}\varphi_{j}(t_{next} )
+\partial_k \partial_i {\nA}_j\varphi_j( t_{next} ), \ \ t<T,\\
\varphi'_{ik}(T)&=\partial_{ik} g(\barX(T)).\\
\endaligned
\end{equation}
Observe that the problem
\eqref{eq:phiend}-\eqref{eq:phiinner}
shows that $\varphi(t_{next})$
depends only on the point values
$$
\{\barX(s): \  t_{next}\le s\le T\},
$$
so that
\begin{equation}\label{eq:2.40}
\partial_{x_k(t)}\varphi_{j}(t_{next} )
=\partial_{x_p(t_{next})}\varphi_{j}(t_{next})\partial_{x_k(t)}\barX_p(t_{next}).
\end{equation}
Finally, the definitions of $\barX$ and $\nA$
in \eqref{eq:2.1} and \eqref{eq:locsolop} imply
\begin{equation}
\label{eq:2.41}
\partial_k {\nA}_p(t,\barX(t))
=\partial_{ x_k(t)}\barX_p(t_{next}),
\end{equation}
which together with (\ref{eq:2.39}-\ref{eq:2.40}) prove that $\varphi'$
satisfies the recursive equation
\begin{equation}\label{eq:phippdyn}
\begin{split}
\varphi'_{ik}(t)&=\partial_i {\nA}_j \ \
\varphi'_{jp}(t_{next})
\partial_k {\nA}_p(t,\barX(t))
+\partial_k \partial_i {\nA}_j \ \ \varphi_j( t_{next} ), \ \ t<T,\\
\varphi'_{ik}(T)&=\partial_{ik} g(\barX(T)),\\
\end{split}
\end{equation}
which is equivalent to \eqref{eq:phiend}-\eqref{eq:phiinner}.
\end{proof}
\begin{remark}\label{rem:meas_trick}
The measurability of
$\bigl(a_{i}(t_{n+1}^{-},\barX(t_{n+1}^{-}))-\bara_{i}(t_n,\barX(t_n))\bigr)\in
{\cF_{t_{n+1}^{-}}}$
proves that for $t=t_{n+1}^{-}$ and any random variable
$\beta$ we have
$$
\aligned
E\bigl[
\bigl(
a_{i}(t,\barX(t))-\bara_{i}(t_n,\barX(t_n))
\bigr)
E[\beta\mid  {\cF_t}]\bigr]
&=E\bigl[E[\bigl(a_{i}(t,\barX(t))-\bara_{i}(t_n;\barX)\bigr)\beta\mid {\cF_t}]\bigr]\\
&=E\bigl[\bigl(a_{i}(t,\barX(t))-\bara_{i}(t_n;\barX)\bigr)\beta\bigr],\\
\endaligned
$$
and in a completely similar way we obtain the same result for the diffusion
terms, i.e.
\begin{equation*}\label{eq:2.36}
E\bigl[\bigl(d_{ij}(t,\barX(t))-\bard_{ij}(t_n;\barX)\bigr)
E[\beta\mid  {\cF_t}]\bigr]
=E\bigl[\bigl(d_{ij}(t,\barX(t))-\bard_{ij}(t_n;\barX))\bigr)\beta\bigr].\\
\end{equation*}
\end{remark}

The proof of Theorem \ref{thm:2.2} is now concluded by combining
Lemmas \ref{lem:2.3}, \ref{lem:2.4},\ref{lem:2.5}, Remark \ref{rem:meas_trick},
and the Central Limit Theorem
to estimate $I_M$ and $II_M$. We have
$$
\aligned
&\int_{\ttt_\ttn}^{\ttt_{\ttn+1}}
E\bigl[\bigl(a_{i}(t,\barX(t^-))-\bara_{i}(t;\barX)
\bigr)
\partial_{i}u(t,\barX(t^-))\bigr] dt\\
=& E\bigl[\sum_{n \in \cJ_\ttn}
          \bigl(
                a_{i}(t_{n+1},\barX(t_{n+1}^{-}))-a_{i}(t_n,\barX(t_n) )
           \bigr)
\partial_{i}u(t_{n+1},\barX(t_{n+1}^{-}))\f{\Delta t_n}{2}\bigr] \\
&+ \sum_{m=0}^{N-1}\OM\bigl((\Delta \ttt_\ttn)^3\bigr)\\
= &
E\bigl[\sum_{n \in \cJ_\ttn}
       \bigl(
             a_{i}(t_{n+1},\barX(t_{n+1}^{-}) )-a_{i}(t_n,\barX(t_n)
       \bigr)
       \partial_{i}\bar u(t_{n+1},\barX(t_{n+1}^{-}))\f{\Delta t_n}{2}
\bigr]\\
& + \sum_{\ttn=0}^{N-1}
 E[\sum_{n \in \cJ_\ttn}(\Delta t_n)^2 \int_{t_{n+1}}^T \OM(\Delta t(s)) ds]
+ \sum_{\ttn=0}^{N-1} \OM\bigl((\Delta \ttt_\ttn)^3\bigr)\\
= &
E\bigl[\sum_{n \in \cJ_\ttn}
       \bigl(
             a_{i}(t_{n+1},\barX(t_{n+1}^{-}) )-a_{i}(t_n,\barX(t_n)
       \bigr)
       E[\varphi_i(t_{n+1-})\mid\cF_{t_{n+1}^{-}}]
      \f{\Delta t_n}{2}
\bigr]\\
& + \sum_{\ttn=0}^{N-1}
 E[\sum_{n \in \cJ_\ttn}(\Delta t_n)^2 \int_{t_{n+1}}^T \OM(\Delta t(s)) ds]
+ \sum_{\ttn=0}^{N-1} \OM\bigl((\Delta \ttt_\ttn)^3\bigr)\\
= &
E\bigl[\sum_{n \in \cJ_\ttn}
       \bigl(
             a_{i}(t_{n+1},\barX(t_{n+1}^{-}) )-a_{i}(t_n,\barX(t_n)
       \bigr)
       \varphi_i(t_{n+1-})\f{\Delta t_n}{2}
\bigr]\\
& + \sum_{\ttn=0}^{N-1}
 E[\sum_{n \in \cJ_\ttn}(\Delta t_n)^2 \int_{t_{n+1}}^T \OM(\Delta t(s)) ds]
+ \sum_{\ttn=0}^{N-1} \OM\bigl((\Delta \ttt_\ttn)^3\bigr)\\
\endaligned
$$
The expansion of the diffusion term appearing in \eqref{eq:rhodef} follows analogously.
\end{proof}
\par
Observe that the number of realizations to determine a reliable error estimate is in
general $\tol^{-1}$, much smaller than the, proportional to $\tol^{-2}$,
number of realizations to approximate $E[g(X(T))]$. For more details on this and the
statistical approximation of the error density $\rho$ see Remark 2.7 in \cite{STZ}.

\section{An Error Estimate with Stochastic Time Steps}\label{Sec:StochTimStep}
This section derives error estimates with time steps
which are stochastic and determined individually
for each realization  by the whole solution
path $\barX$. The analysis will
use the Malliavin derivative, $\partial_{W(t)}Y$, which is
the first variation of a process $Y$ with respect to a perturbation
$dW(t)$,  at time $t$ of the Wiener process, cf. \cite{Nu}.
The Malliavin derivative
for a stochastic integral $X$ is related to the first variation,
$\partial_{x(t)}\barX$, for a perturbation
of the  position at time $t$ by
\begin{equation}\label{eq:3.1}
\begin{split}
\partial_{W^\ell(t)}X(\tau)&=\tfrac{\partial X_k(t)}{ \partial W^\ell(t)}
\partial_{x_k(t)} X(\tau) = b_k^\ell(X(t))\partial_{x_k(t)}X(\tau), \quad\tau>t,\\
\partial_{W^\ell(t)}X(\tau)&=0 \quad\tau<t,\\
\end{split}
\end{equation}
if $dX_k= a_k(X(t)) dt + b_k^\ell(X(t)) dW^\ell$ (cf. \eqref{eq:2.15}).
\par
We shall restrict the analysis to time steps which are constructed
by first sampling the jump times to augment an a priori given
time-discretization 
$\Delta \tilde{t}$, 
obtaining 
$\Delta t[0]$ 
(see \eqref{eq:zero}) and then using the refinement
criterion
\begin{equation}
\aligned
\Delta t(t) &= \Delta t[0](t)/2^n, \ \ \hbox{for some natural number} \ n=n(t,\omega),\\
|\rho(t,\omega)|\bigl(\Delta t(t)\bigr)^2 &< \hbox{constant},\\
\endaligned\label{eq:3.2}
\end{equation}
with an approximate error density function, $\rho$, satisfying,
 for $s\in [0,T],
 \ t\in [0,T]$ and all outcomes $\omega$,
 the
uniform upper and lower bounds
\begin{equation}\label{eq:3.3}
\aligned
c(\tol)\le |\rho(s,\omega)|&\le C(\tol),\\
|\partial_{W(t)}\rho(s,\omega)|&\le C(\tol),\\
\endaligned
\end{equation}
for some positive functions $c$ and $C$, with $\tol/c(\tol)\rightarrow 0$
as $\tol\rightarrow 0$.
For each realization
successive subdivisions of the steps
yield the largest time steps satisfying \eqref{eq:3.2}.
The corresponding stochastic increments $\Delta W$
will have the correct distribution,
with the necessary independence,
if the increments $\Delta W$ related to the new steps are
generated by Brownian bridges, cf. \cite{KS},
i.e. the time steps are generated by
conditional expected values of the Wiener process.

Let $\delta$ be a constant
approximating 
$\frac{\tol}{E[N]}$, 
where 
$E[N]$ 
is the expected
number of steps.
The analysis in this section with
adaptive non adapted time steps, satisfying \eqref{eq:3.2}-\eqref{eq:3.3},
is based on the following Stochastic time step algorithm
described in next page.

Section \ref{sec:sto} presents a more precise formulation of this
algorithm.
Lemma \ref{lem:3.1} and Theorem \ref{thm:stoch_err_exp}
below show that although the steps
generated by \eqref{eq:3.2}-\eqref{eq:3.3} through the algorithm above
are not adapted, the method indeed converges to the correct limit as the
forward Euler method with adapted time steps.
\begin{lemma}[Strong convergence]\label{lem:3.1}
Suppose that $a,b,g, X$ satisfy the assumptions in Lemma \ref{lem:2.1}
and that $\barX$ is constructed by the forward Euler method,
based on  the stochastic time step algorithm above,
with step sizes $\Delta t_n$ satisfying (\ref{eq:3.2}-\ref{eq:3.3})
and their
corresponding $\Delta W_n$ are generated by Brownian bridges.
Assume also that $\barX(0)=X(0)$. Then
$$
\sup_{0\le t\le T}\sqrt{E[|X(t)-\barX(t)|^2]}
= \OM(\sqrt{\Delta t_{\hbox{sup}}})
= \OM(\sqrt{\frac{\tol}{c(\tol)}})\rightarrow 0,
$$
as $\tol\rightarrow 0$, where $\Delta t_{\hbox{sup}}\equiv
\sup_{n,\omega}\Delta t_n(\omega)$.
\end{lemma}
\begin{proof}
Let $\sigmatt$
be the $\sigma$-algebra generated by the jumps and marks constructed in Remark \ref{Remarkk_1}.
Consider the conditional expectation $E[|X(t)-\barX(t)|^2\mid\sigmatt]$
and apply Lemma 3.1 from \cite{STZ}, using also that there is no time
discretization error at the jump nodes.
\end{proof}
\par
In addition to the dual functions $\varphi$ and $\varphi'$
in Theorem~\ref{thm:2.2},
the new error expansion for stochastic time steps in
Theorem~\ref{thm:stoch_err_exp}
below also uses, for $t = t_{n+1}^{-}$, the discrete dual variation
\begin{equation}\label{eq:3.9}
\varphi''_{ikm}(t)\equiv \partial_{x_m(t)}\varphi'_{ik}(t)
\equiv\tfrac{\partial\varphi'_{ik}(t;\barX(t)=x)}{\partial x_m},
\end{equation}
which satisfies the backward problem
\eqref{eq:phiend} -- \eqref{eq:phiinner}, i.e.,
\begin{lemma}\label{lem:third_var}
Let $\nA$, $\varphi$ and $\varphi'$
be defined by  \eqref{eq:locsolop}
and \eqref{eq:phiend} -- \eqref{eq:phiinner}.
Then $\varphi''$  is given by
\begin{equation}\label{eq:3.10}
\aligned \varphi''_{ikm}(t)&=\partial_i \nA_j(t,\barX(t))\partial_k \nA_p(t,\barX(t))
\partial_m \nA_r(t,\barX(t))
\varphi''_{jpr}(t_{next} )\\
&\quad +\partial_{im} \nA_j(t,\barX(t))\partial_k \nA_p(t,\barX(t))
\varphi'_{jp}(t_{next} )\\
&\quad +\partial_i \nA_j(t,\barX(t))\partial_{km}\nA_p(t,\barX(t))
\varphi'_{jp}(t_{next} )\\
&\quad +\partial_{ik} \nA_j(t,\barX(t))\partial_m \nA_p(t,\barX(t))
\varphi'_{jp}( t_{next} )\\
&\quad +\partial_{ikm}\nA_j(t,\barX(t))
\varphi_{j}( t_{next} ), \ \ t<T,\\
\varphi''_{ikm}(T)&=\partial_{ikm} g(\barX(T)).\\
\endaligned
\end{equation}
which is equivalent to \eqref{eq:phiend} -- \eqref{eq:phiinner}
with $t_{next}$ as in Definition~\ref{def:ldso}.
\end{lemma}
\noindent

\fbox{
\begin{minipage}[h]{12.8cm}
{\tt Stochastic time step algorithm:}

{\bf Do } for $M$ realizations $\omega_j, \ j=1,\ldots,M$:

\hspace{0.5cm} Sample jump times:
$0<\tau_1(\omega_j)<\cdots<\tau_{\widehat{N}(\omega_j)}(\omega_j)< T$.

\hspace{0.5cm}
Consider an augmented partition given by the union
\begin{equation}\label{eq:zero}
\Delta t[0]=\{t_n(\omega_j)\}_{n=0}^{\ssy N_{\ssy A}(\omega_j)}=\{\ttt_n\}_{n=0}^{\ssy N}
\cup\{\tau_n(\omega_j)\}_{n=0}^{\ssy \widehat{N}(\omega_j)}
\end{equation}

\hspace{0.5cm} that has  $N_{\ssy A}(\omega_j)= N+\widehat{N}(\omega_j)$ ($a.s.$) time steps.

\noindent
\hskip0.5truecm{\vbox{\hsize12.0truecm\noindent\strut
{\tt STEP 1}: Set $k=0$. Start with the initial coarse mesh $\Delta t[0]$
and compute $\Delta W[0]$.
\strut}}

\noindent
\hskip0.5truecm{\vbox{\hsize12.0truecm\noindent\strut
{\tt STEP 2}: { For} the piecewise constant mesh function
$\Delta t[k]$ with corresponding noise $\Delta W[k]$,
compute $\barX[k]$ and
the weight function $\rho[k]$ defined in (\ref{eq:3.12a}-\ref{eq:3.12b}).
\strut}}

\noindent
\hskip0.5truecm{\vbox{\hsize12.0truecm\noindent\strut
{\tt STEP 3}: Define $r(t)\equiv  |\rho[k](t)|(\Delta t[k](t))^2$
and let for all $t$
$$
\Delta t[k+1](t)= \left\{
\aligned
&\Delta t[k](t),\  \ \ \ \ \hbox{if}\  r(t)<\delta,  \ \ (\dagger)\\
&\Delta t[k](t)/2, \  \ \hbox{if}\  r(t)\ge \delta,
\ \ (\star)\\
\endaligned
\right.
$$
and in the refinement case ($\star$) construct $\Delta W[k+1]$
by Brownian
bridges based on the already known $\Delta W[k]$.
\strut}}
%

\noindent
\hskip0.5truecm{\vbox{\hsize12.0truecm\noindent\strut
{\tt STEP 4}:
{\bf If} at least one step of $\Delta t[k]$ is refined by ($\star$), increment
$k$ by $1$ and {\bf goto} { Step 2}.
{\bf Else} all steps of $\Delta t[k]$
satisfy ($\dagger$) and accept the approximation  $g(\barX(T,\omega_j))$
and {\bf goto} the next realization of $p$ and $W$.
\strut}}
\break
{\bf Enddo }

{\vbox{\hsize12.0truecm\noindent\strut
{\bf If} the statistical error,
$E[g(\barX(T))]-
\frac{1}{M}\,\sum_{j=1}^M g(\barX(T,\omega_j))$,
is sufficiently small {\bf stop},
{\bf else} restart with a larger $M$.
{\bf Endif}
\strut}}
\end{minipage}}
\vskip2mm
\begin{proof}
Differentiation of the backward recursive equation
\eqref{eq:phippdyn}
and the relations (\ref{eq:2.37}-\ref{eq:2.41})
together with
\begin{equation}
\partial_{x_m(t)}\varphi'_{jp}(t_{next} )=\partial_{x_r(t_{next})}\varphi'_{jp}(t_{next} )
\partial_{x_m(t)}\barX_r(t_{next}),
\label{eq:3.11}
\end{equation}
prove as in \eqref{eq:2.38}-\eqref{eq:2.41}
that $\varphi''$ satisfies \eqref{eq:3.10}.
Here, \eqref{eq:3.11} holds since the linear system
for the variable $\bigl(\varphi(t_{next}),\varphi'(t_{next})\bigr)$
depends only on the point values $\{\barX(s): \  t_{next}\le s \le T\}.$
\end{proof}
\par
The following theorem derives an error estimate
applicable both to adaptive deterministic time steps
and to the stochastic time step algorithm; the
assumptions and the proof of the theorem focus
on stochastic steps, however a modification
to deterministic time steps is straightforward.
The computable
error density $|\tilde\rho|$ of this error estimate
can then be cut-off for small and large values
to satisfy \eqref{eq:3.3}, see \eqref{eq:density} and \eqref{eq:density_D}.


\begin{thm}[Stochastic time steps error expansion]
\label{thm:stoch_err_exp}
Suppose that $a,b,g, X$ satisfy the assumptions in Lemma \ref{lem:2.1}
and that $\barX$ is constructed by the forward Euler method
with step sizes $\Delta t_n$ satisfying (\ref{eq:3.2}-\ref{eq:3.3}) and the
corresponding $\Delta W_n$ are generated by Brownian bridges,
following the stochastic time step algorithm in Lemma \ref{lem:3.1}.
Assume also that $\barX(0)=X(0)$ and $E[|X(0)|^{k_0} ]\le C$
for some $k_0\ge 16$.
Then the time discretization error has the following expansion,
 based on both the drift and diffusion fluxes and the
discrete dual functions $\varphi,\ \varphi'$ and $\varphi''$ given in
\eqref{eq:phiend} -- \eqref{eq:phiinner},
with computable leading order terms
%
%
\begin{equation}\label{eq:3.12a}
\aligned E[g(X(T))-g(\barX(T))]
&=E\Bigl[\sum_{n=0}^{N_{\ssy A}-1}\tilde\rho(t_n,\barX)(\Delta t_n)^2\Bigr]\\
&\quad + \OM\Bigl( \sqrt{\frac{\tol}{c(\tol)}}
\Bigl(\frac{C(\tol)}{
c(\tol)}\Bigr)^{8/k_0}\Bigr)
E\Bigl[\sum_{n=0}^{N_{\ssy A}-1}(\Delta t_n)^2 \Bigr],\\
\endaligned
\end{equation}
where
\begin{equation}\label{eq:3.12b}
\aligned \tilde\rho(t_n,\barX)&\equiv {\half}\Bigl( \bigl(
\pa_t a_k+\partial_j a_ka_j +\partial_{ij}a_k
d_{ij}
\bigr)\varphi_k(t_{n+1}^{-})\\
&
\quad+\bigl( \pa_t d_{km}+\partial_j d_{km}a_j
+ \partial_{ij}d_{km} d_{ij}
+2\partial_ja_k d_{jm}
\bigr)\varphi'_{km}(t_{n+1}^{-})\\
&\quad
+\bigl( 2\partial_j d_{km}d_{jr}
\bigr)\varphi''_{kmr}(t_{n+1}^{-})\Bigr),\\
\endaligned
\end{equation}
and the terms in the sum of \eqref{eq:3.12b} are evaluated at the
a posteriori known points $(t_n,\barX(t_n))$, i.e.
$$\aligned
&\partial_\alpha a \equiv \partial_\alpha a(t_n,\barX(t_n)),\\
&\partial_\alpha b \equiv \partial_\alpha b(t_n,\barX(t_n)),\\
&\partial_\alpha d \equiv \partial_\alpha d(t_n,\barX(t_n)).\\
\endaligned
$$
\end{thm}
\begin{proof}
%
%
We consider the difference
\begin{equation*}
g(X(t))-g(\barX(t))
\end{equation*}
and apply Theorem~3.3 from \cite{STZ}, using also that there is
no time discretization error at the jump nodes.
To this end, denote the set of stochastic time nodes by
$\cJ \equiv \{0= t_0,t_{1}^{-}, t_{1},t_{2}^{-}, t_2,\ldots,t_N =T\}$ and
recall that the notation 
\begin{equation}\label{eq:3.13a}
\barX(t_{next})=\nA(\barX(t)),\quad t\in\cJ,
\end{equation}
introduced in Definition~\ref{def:ldso} denotes one step with the Euler method.
Write
similarly one step with the exact solution
\begin{equation}
\label{eq:3.13b}
X(t_{next})=\hat{\nA}(X(t)),\ \ t\in \cJ.
\end{equation}
Introduce the notation $X^t\equiv X(t)$
and $\barX^t\equiv \barX(t)$.
Now verify the representation
\begin{equation}
\label{eq:3.14}
g(X(T))-g(\barX(T)) = \sum_{n=0}^{N_{\ssy A}-1} (\hat{\nA}(\barX(t_n))-\nA(\barX(t_n)))_i\tilde\varphi_i(t_{n+1}^{-})
\end{equation}
where the weight functions are defined recursively by the linear backward recursion
\begin{equation}
\label{eq:3.15}
\aligned
&\tilde\varphi_i(T)
=\int_0^1\partial_ig(sX(T)+(1-s)\barX(T))ds,\\
&\tilde\varphi_i(t)
=\Bigl(
\int_0^1\partial_i
\hat{\nA}_j(sX(t)+(1-s)\barX(t))ds\Bigr)\tilde\varphi_j(t_{next}),
\
t \in \cJ
.\\
\endaligned
\end{equation}
To verify \eqref{eq:3.14}, first observe that by construction
of the Euler method at every jump point $t = t_{n}^{-}$ there is no local error
in computing the next $\barX$ value at time $ t_n$, i.e.
$$
\hat{\nA}(\barX^t)-\nA(\barX^t) =0
$$
so \eqref{eq:3.14} is equivalent to
\begin{equation}\label{eq:3.14c}
g(X(T))-g(\barX(T))= \sum_{t\in \cJ} (\hat{\nA}(\barX^t)-\nA(\barX^t))_i\tilde\varphi_i(t_{next}).
\end{equation}
Then telescoping cancelation gives
\begin{equation}\label{eq:3.16}
g(X(T))-g(\barX(T))=
\sum_{t\in \cJ}
\bigl( (X^{t_{next}}-\barX^{t_{next}})_i\tilde\varphi_i(t_{next})
- (X^{t}-\barX^{t})_i\tilde\varphi_i(t)\bigr).
\end{equation}
Use the definitions (\ref{eq:3.13a},\ref{eq:3.13b}) and split the first term in the
sum of \eqref{eq:3.16} into
$$
(\hat\nA(X^{t})-\hat\nA(\barX^t) + \hat\nA(\barX^t)
- \nA(\barX^t))_i\tilde\varphi_i(t_{next}).
$$
The two first terms above and the last term in the sum of
\eqref{eq:3.16} combine to zero
by \eqref{eq:3.15}:
\begin{equation*}
(\hat\nA(X^{t})-\hat\nA(\barX^t))_i\tilde\varphi_i(t_{next})(X^{t}-\barX^{t})_i\tilde\varphi_i(t),
\end{equation*}
which proves \eqref{eq:3.14}.
\par
The next step is to use the Malliavin derivative to
analyze the  expectation of the representation
\eqref{eq:3.14} by studying the dependence of $\barX$ and
$\tilde\varphi$ on a small increment $dW$. This follows exactly
the lines of the proof of Theorem~3.3 from \cite{STZ},
and it is not reproduced here.
\end{proof}
%
%
%
\section{Adaptive time-stepping algorithms}\label{sec:algo}
%
Here, we describe two adaptive time-stepping algorithms for the
weak approximation problem of \eqref{eq:1.1} based on the
approximation error formulas described in the previous section.
They are very similar to those introduced in \cite{MSTZ}.
Algorithm-D is based on a quasi-deterministic mesh that is fixed for all
realizations and its adaptive strategy is based on averaged
information from the a posteriori error formula. On the other
hand, Algorithm-S can adapt the time discretization differently
for each realization. Proper sample of the Wiener process is
possible by means of Brownian bridges. Both adaptive algorithms
choose adaptively the number of realizations and the size of time
steps to efficiently bound the approximation error by a prescribed
error tolerance.
\subsection{Computational error splitting }
The weak approximation computational error
of the Monte Carlo Euler method,
$\errorc\equiv E\big[g(X(T))\big]
-\tfrac{1}{M}\sum_{j=1}^{M}g\big(\barX(T,\omega_j)\big)$,
naturally separates to the time discretization error
$\errort\equiv E\big[g(X(T))\big]
-E\big[g(\barX(T))\big]$
and the statistical error
$\errors\equiv E\big[g(\barX(T))\big]
-\tfrac{1}{M}\sum_{j=1}^{M}g\big(\barX(T,\omega_j)\big)$, i.e.,
$\errorc=\errort+\errors$.
Thus, the control the computational error is related to the
a combined control of the time discretization error via the
choice of the time steps $\Delta{t}$ and of the statistical
error via the choice of the number $M$ of the realizations.
Therefore, we split a given computational error tolerance,
$\tol>0$, into a statistical error tolerance $\tol_{\ssy S}$ and a
time discretization error tolerance $\tol_{\ssy T}$ (see
\cite{STZ}, \cite{MSTZ}) by
\begin{equation}\label{eq:tol}
\tol_{\ssy T}=\tfrac{1}{3}\,\tol\quad\mbox{\rm and}
\quad\tol_{\ssy S}=\tfrac{2}{3}\tol.
\end{equation}
%
%
%
\subsection{Control of the statistical error}
For $M$ independent samples $\{ Y(\omega_j) \}_{j=1}^{\ssy M}$ of a
random variable $Y$, with $E\big[\,|Y|^6\,\big]<\infty$, 
define the sample average by 
%
$${\ma}(Y;M) \equiv \frac{1}{M} {\sum_{j=1}^M} Y(\omega_j)$$
and the sample standard deviation by
$$
{\mathcal S}(Y;M) \equiv
\left\{{\ma}(Y^2;M)-({\ma}(Y;M))^2\right\}^{\frac{1}{2}}.
$$
Let $\sigma_{\ssy Y}\equiv
\left\{E\big[\big|Y-E[Y]\big|^2\big]
\right\}^{\frac{1}{2}}$
and
$Z_{\ssy M}\equiv\frac{\sqrt{M}}{\sigma_Y}\left(\ma(Y;M)-E[Y] \right)$
with cumulative distribution function $F_{Z_M}(x) \equiv P(Z_M
\leq x), ~x\in \rset$. Let $\lambda\equiv
(E[|Y-E[Y]|^3])^{1/3}/\sigma_Y<\infty,$ then the Berry-Esseen
theorem, cf. \cite{Du}, gives the following estimation in the
Central Limit Theorem
\[
\sup_{x\in\rset}|F_{Z_M}(x)-\normal(x)| \le
\frac{3}{\sqrt{M}}\lambda^3
\]
for the rate of convergence of $F_{Z_M}$ to the distribution
function $\normal$ of a normal random variable with mean zero and
variance one, i.e.
\begin{eqnarray}
\normal(x) = \int_{-\infty}^x \frac{1}{\sqrt{2 \pi}} \exp \left( -
\frac{1}{2} s^2\right) ds. \label{eq:normal}
\end{eqnarray}
Since in the examples below $M$ is sufficiently large, i.e. $M\gg
36 \lambda^6$, the statistical error
$$\errors(Y;M)\equiv E[Y]-\ma(Y;M)$$ satisfies,
by the Berry-Esseen theorem, the following probability
approximation
$$P\left( \left[ |\errors(Y;M)|
\leq{c_{_0}}\frac{\sigma_Y}{\sqrt{M}}\right]\right) \simeq
2\normal(c_{_0})-1.$$
In practice choose some constant $c_{_0}\ge1.65$, so the normal
distribution satisfies $1>2\normal(c_{_0})-1\ge 0.901$ and the event
\begin{eqnarray}
| \errors(Y;M) |\leq \mathtt{E}_S (Y;M) \equiv
c_{_0}~\frac{{\mathcal{S}}(Y;M)}{\sqrt{M}} \label{eq:static_bdd}
\end{eqnarray}
has probability close to one, which involves the additional step
to approximate $\sigma_Y$ by $\mathcal{S}(Y;M)$, cf. \cite{Fi}.
Thus, in the computations $\mathtt{E}_S(Y;M)$ is a good
approximation of the statistical error $\errors(Y;M)$.
\par
For a given $\tol_{\ssy S}>0$, the goal is to find $M$ such that
$\ErrS(Y;M)\leq\tol_{\ssy S}$. The following algorithm adaptively
finds the number of realizations $M$ to compute the sample average
${\ma}(Y;M)$ as an approximation to $E[Y]$. With probability
close to one, depending on $c_{_0}$, the statistical error in the
approximation is then bounded by $\tol_{\ssy S}$. For technical
reasons (see \cite{MSTZ}) we choose $M=2^n$, $n\in\nset$.
%
\par\medskip

{\tt routine Monte-Carlo}($\tol_S$, $M_0$; $EY$)

\hspace{1cm}Set the batch counter $m=1$, $M[m]=M_0$ and
$\ErrS[m]=+\infty$.

\hspace{1cm}\textbf{Do while} ($\ErrS[m]>\tol_S$)
\begin{equation}
\label{MC}\begin{array}{l} \mbox{Compute $M[m]$
$\underline{\text{\rm new}}$ samples of $Y$, along
        with the sample}\\

\mbox{average $EY \equiv {\ma}(Y;M[m])$, the sample standard
deviation
        }\\

\mbox{${\mathcal{S}}[m]\equiv{\mathcal{S}}(Y;M[m])$ and the statistical error estimation }\\

\mbox{$\ErrS[m+1]\equiv\ErrS(Y,M[m])$. Compute $M[m+1]$ by }\\
\mbox{{\tt change$\underline{\kenn}M$}
        ($M[m]$, ${\mathcal{S}}[m]$, $\tol_S$; $M[m+1]$).}\\

\mbox{Increase $m$ by 1.}
    \end{array}
\end{equation}
\hspace{1cm}\textbf{end-do}

{\tt end of Monte-Carlo}

\medskip

{\tt routine change$\underline{\kenn}M$} ($M_{in}$,
${\mathcal{S}}_{in}$, $\tol_S$; $M_{out}$)
\begin{equation}
\label{eq:choose_M}
    \begin{array}{ll}
    M^* &=\min\Bigl\{\text{\rm integer part}\kenn
    \Bigl(\frac{c_{_0}\kenn{\mathcal{S}}_{in}}
    {\kenn\tol_S}\Bigr)^2,
    \kenn\mbox{MCH}\times M_{in}\Bigr\}   \\
    n & = \text{integer part}~ (\log_2 M^*) +1 \\
    M_{out}& = 2^n.
    \end{array}
\end{equation}

{\tt end of change$\underline{\kenn}M$}

\medskip

Here, $M_0$ is a given initial value for $M$, and $\mbox{MCH}>1$
is a positive integer parameter introduced to avoid a large new
number of realizations in the next batch due to a possibly
inaccurate sample standard deviation $\mathcal{S}[m]$. Indeed,
$M[m+1]$ cannot be greater than $\mbox{MCH}\times{M[m]}$.
%
%
%
%
\subsection{Deterministic time stepping algorithm}\label{sec:determ}
Following closely \cite{MSTZ} and \cite{MSST}, we present an
adaptive algorithm based on cut-off of the error density
${\widetilde\rho}$ in \eqref{eq:3.12b} of
Theorem~\ref{thm:stoch_err_exp},
$\rho_{\ssy D}$,  which is defined
as

\begin{equation}\label{eq:density_D}
\rho_{\ssy D}^n\equiv \min\left( \max\left(
\Big|\tfrac{1}{(\Delta{\ttt}_n)^2}
\,\sum_{\ell\in{\mathcal J}_n}
(\Delta{t}_{\ell})^2
{\widetilde\rho}(t_{\ell},\barX)\Big|,
\tol^{\frac{1}{9}}\right),\tol^{-1}\right), \quad n=1,\dots, N.
\end{equation}
The error expansion in Theorem~~\ref{thm:stoch_err_exp} motivates us to
approximate
the time discretization error  by
\begin{equation}\label{eq:approx_error}
|\errort|\lesssim{E}\left[
\sum_{n=1}^{N} r_n\right]
\end{equation}
where the error indicator, $r_n$, is defined by
\begin{equation}\label{eq:indicator}
r_n\equiv
\rho^n_{\ssy D}\,(\Delta\ttt_n)^2
\quad n=1,\dots,N.
\end{equation}

The main advantage of the deterministic time stepping algorithm over
the stochastic time stepping algorithm is that the number $M_T$ of realizations
necessary to determine the optimal deterministic stepping scheme
is considerable smaller than $M$, whereas, 
in the stochastic time stepping algorithm, 
a refinement of the partition is carried out in each one of the $M$ trajectories, leading to 
a considerable larger amount of computational work.


As pointed out before, the error expansion derived in Theorem \ref{thm:stoch_err_exp}
is also valid with deterministic  time steps. Therefore, we have some flexibility in the  choice
of error densities for the adaptive algorithm with deterministic time steps
because we can also use the results from Theorem  \ref{thm:2.2}.
There are some practical differences to mention. On the one hand, the variance
of the averaged error density from \eqref{eq:rhodef} is $\mathcal{O}(\frac{1}{\dtt_n M_T})$.
This feature has been observed in \cite{MSTZ2} and \cite{STZ}, where a local filtering procedure
was proposed to reduce the variance of the error density estimator. A positive
feature of this error density is that it does not require the computation of the second variation, $\varphi''$, which may be computationally expensive for large $d$.
On the other hand, the averaged error density from \eqref{eq:3.12b} has a much smaller
variance $\mathcal{O}(\frac{1}{M_T})$ which does not need filtering but it requires
the computation of $\varphi''$.
In this work we will not discuss further this choice and only show numerical results with adaptive deterministic time steps
based on the error density \eqref{eq:3.12b}.

The approximation of the time discretization error in the
right hand side of \eqref{eq:approx_error} can be separated into
two parts
\begin{eqnarray}
\ken\ken\ken E \left[ \sum_{n=1}^{N} r_n\right] \leq \ma\left(
\sum_{n=1}^{N} r_n ;M_T\right) + \left|E \left[ \sum_{n=1}^{N}
r_n\right] -
    \ma\left( \sum_{n=1}^{N} r_n ;M_T\right)\right|,
\label{eq:time_split}
\end{eqnarray}
where the second error term in the right hand side of
\eqref{eq:time_split} is with probability close to one
asymptotically bounded by
\begin{eqnarray}
\left|E \left[ \sum_{n=1}^{N} r_n\right] -
    \ma\left(\sum_{n=1}^{N} r_n; M_T \right)\right|
\lesssim  ~\mathtt{E}_{TS}~ \equiv c_0
    \frac{\mathcal{S}\left( \sum_{n=1}^{N} r_n; M_T \right)}{\sqrt{M_T}}
\label{eq:E_TS}
\end{eqnarray}
and the first term defines $\mathtt{E}_{TT} \equiv \ma\left(
\sum_{n=1}^{N} r_n ; M_T\right)$. Then for a given $\tol_T>0$, the
goal is to construct a partition $\Delta \ttt$ of $[0,T]$, with as
few time steps and realizations $M_T$ as possible, such that
$\mathtt{E}_{\ssy TT} + \mathtt{E}_{\ssy TS} \leq \tol_{\ssy T}$.
To this end, first split the time discretization tolerance
$\tol_{\ssy T}$ in two positive parts $\tol_{\ssy TT}$ and
$\tol_{\ssy TS}$ for $\mathtt{E}_{\ssy TT}$ and $\mathtt{E}_{\ssy
TS}$, respectively.
The statistical error of the time discretization using the density
\eqref{eq:density_D} is $\mathcal{O}(\frac{\supdt}{\sqrt{M_T}})$.
Therefore the percentage of the tolerance, $\tol$, devoted to the
control of the statistical time discretization error can be
arbitrary small as $\supdt \ra 0$. In practice we choose
\begin{equation}
\tol_{TT} = \frac{2}{3}\tol_T = \frac{2}{9}\tol,
\qquad
\tol_{TS} = \frac{1}{3} \tol_T =\frac{1}{9}\tol. \label{eq:time_tol}
\end{equation}
The control of the statistical time discretization error
determines the number of realizations $M_T$ necessary to ensure a reliable
choice of the time discretization in the deterministic time
stepping algorithm.
Similarly as in \cite{MSTZ}, here it is  optimal to equidistribute
the error contributions from different time intervals. Thus, the
goal of the adaptive algorithm described below is to construct a
deterministic time partition $\Delta \ttt$ of $[0,T]$ such that
\begin{eqnarray}
\bar r_n \equiv {\mathcal{A}(r_n;M_T)}
\leq d_1 \frac{\tol_{TT}} {N}, ~~ ~n=1,\dots,N,
\label{eq:condition_det}
\end{eqnarray}
where $d_1= 2$, see Remark 3.9
in \cite{MSTZ}.

To achieve \eqref{eq:condition_det}, start with an initial partition $\Delta \ttt[1]$
and then specify iteratively a new partition $\Delta{\ttt}[k+1]$,
from $\Delta \ttt[k]$, using the following refinement strategy:

\begin{eqnarray}
&&\mbox{\bf for} ~~n=1,2,\ldots ,N[k]\\
&~~&\mbox{\bf if}~  \bar r_n[k] \geq d_1 \frac{\tol_{TT}}{N[k]},
~\mbox{{ \bf then} divide $\Delta \ttt_n[k]$ into $2$ uniform
substeps.}
\label{eq:divide_det} \\
&~~& \mbox{{\bf else} let the new step be the same as the old}
\label{eq:remain_det} \nonumber \\
&~~& \mbox{\bf endif}.\\
&& \mbox{\bf endfor}. \nonumber
\end{eqnarray}
until the following stopping criteria is satisfied:
\begin{eqnarray}
\label{eq:stop-cond1_det} &&\mbox{\bf if}~ \left( \max _{1\leq n
\leq N[k]} \bar r_n[k] < D_1 \frac{\tol_{TT}}{N[k]} \right)
~~\mbox{{\bf then} stop.}
\end{eqnarray}
Here $D_1$ is a given constant satisfying $D_1>\frac{2}{c}d_1$ where $c \approx 1/2$, see \cite{MSTZ}. The
combination of \eqref{eq:time_split} and \eqref{eq:stop-cond1_det}
asymptotically guarantees a given level of accuracy,
$\mathtt{E}_{TT} < D_1 \tol_{TT}$. The positive numbers $D_1$ is
motivated to avoid slow convergence in case almost all $\bar r_n$
satisfy \eqref{eq:stop-cond1_det}, as in Section~\ref{sec:sto}.
%
\par
Now we are ready for the detailed definition of the adaptive
algorithm with deterministic steps:

\medskip

{\tt Algorithm D}

\textbf{Initialization} Choose:
   \begin{enumerate}
   \item an  error tolerance, $\tol \equiv \tol_S + \tol_{TT} + \tol_{TS}$,
   \item a number, $N[1]$, of initial uniform steps $\Delta \ttt[1]$ for $[0,T]$,
   \item a number, $M[1]$, of initial realizations and set $M_T[1] = M[1]$,
   \item a number, $d_1 = 2$ in \eqref{eq:divide_det} and  $c=1/2+1/20$   to compute $D_1$
    using $D_1>\frac{2}{c}d_1$, and
   \item a constant $c_0 \geq 1.65$ and an integer MCH$\ge 2$ to
    determine the number of realizations in \eqref{eq:choose_M}.
   \end{enumerate}

Set the iteration counter, $k$, for time refinement levels, to 1
and set the statistical

error, ${\mathtt E}_{TS} = + \infty$ and $\bar r [k] = + \infty$.

%
%
\textbf{Do while} ( $\bar r [k]$ violates the stopping
\eqref{eq:stop-cond1_det} or $\mathtt{E}_{TS} > \tol_{TS}$ )

\hspace{0.5cm}Compute the sample averages and the error estimates
on $\dtt[k]$

\hspace{0.5cm}by calling {\tt Euler}. Set $M_T[k+1] = M_T[k]$ and
$\dtt[k+1] = \dtt[k]$.

\hspace{0.5cm}\textbf{If} ( $\bar r[k]$ violates the stopping
            \eqref{eq:stop-cond1_det} )

\hspace{1cm}For all time steps $i = 1, \dots, N[k]$,
          do the refinement process \eqref{eq:divide_det}

\hspace{1cm} to update
          $\Delta \ttt[k+1]$ from $\Delta \ttt[k]$.

\hspace{0.5cm}\textbf{elseif} ( $\mathtt{E}_{TS} > \tol_{TS}$ )

\hspace{1cm}Update $M_T[k+1]$ by
    {\tt change$\underline{\kenn}$M}
    $(M_T[k],{\mathcal{S}}_{TS}[k],\tol_{TS};M_T[k+1])$.

\hspace{0.5cm}\textbf{end-if}

\hspace{0.5cm}Increase $k$ by 1.

\textbf{end-do}

%
%
Compute an approximation, $Eg$, for $E[g(\barX(T))]$ with fixed
time mesh $\Delta \ttt=\Delta \ttt[k]$

by {\tt  Monte-Carlo}($\tol_S$, $M_T[k]$; $Eg$) in \eqref{MC}.

Accept $Eg$ as an approximation of $E[g({{X}}(T))]$, since the
estimate of the

computational error is bounded by $\tol$.

\bigskip

{\tt routine Euler}

\hspace{0.5cm}For each  $M_{_T}[k]$ new realizations, sample jump
times with their corresponding marks.

\hspace{0.5cm}As described in Section \ref{Kapote}, use $\Delta{\ttt}[k]$ to
compute corresponding realizations of the Euler method.

\hspace{0.5cm} Update the approximations of the time
discretization error

\hspace{0.5cm}indicators $\bar r[k]$ and the statistical time
discretization error $\ETS[k]$ and

\hspace{0.5cm}compute the sample standard deviation
$\mathcal{S}_{TS}[k] \equiv \mathcal{S}(g(\barX(T)); M_T[k])$.

{\tt end-of-Euler}

\subsection{The stochastic time stepping algorithm}\label{sec:sto}
Now we describe an adaptive algorithm with stochastic time steps
 based on a cut-off of the error density
${\widetilde\rho}$ introduced in \eqref{eq:3.12b} of
Theorem~\ref{thm:stoch_err_exp},
$\rho_{\ssy S}$, defined by
as
\begin{equation}\label{eq:density}
\rho_{\ssy S}^n
\equiv \min\left(\max\left(|{\widetilde\rho}(t_n,\barX)|,\tol^{\frac{1}{9}}\right),\tol^{-1}\right),
\quad n=1,\dots, N_{\ssy A}.
\end{equation}
Following the error expansion in Theorem~~\ref{thm:stoch_err_exp},
the time discretization error is approximated by
\begin{equation}\label{eq:approx_errors}
|\errort|\lesssim{E}\left[
\sum_{n=1}^{N_{\ssy A}} r_n\right]
\end{equation}
where the error indicator, $r_n$, is defined by
\begin{equation}\label{eq:indicators}
r_n\equiv
\rho^n_{\ssy S}\,(\Delta t_n)^2
,\quad n=1,\dots,N_{\ssy A}.
\end{equation}
In this case it is optimal, cf. \cite{MSTZ}, to equidistribute the
error contributions among all time steps and all realizations.
In other words,  the goal of the adaptive algorithm is to construct a
time partition $\Delta t$ of $[0,T]$ for each realization such
that
\begin{eqnarray}
r_n \leq s_1 \tfrac{\tol_{T}} {E[N_{\ssy A}]}, ~~~
~n=1,\dots,N_{\ssy A}, \label{eq:condition}
\end{eqnarray}
where $s_1 =2$, see Remark 3.1 in \cite{MSTZ}.  Note that in practice the quantity $E[N_{\ssy A}]$ is
not known and we can only estimate it by a sample average
$\ma(N_{\ssy A};M)$ from the previous batch of realizations.
The statistical error $|E[N_{\ssy A}] - \ma(N_{\ssy A};M)|$ is then bounded by
$\mathtt{E}_S(N_{\ssy A};M)$, with probability close to one, by the same
argument as in \eqref{eq:static_bdd}.

Let $\barN_{\ssy A}[j]\equiv\ma(N_{\ssy A};M[j])$ be the sample average of the
final number of time steps in the $j$-th batch of $M[j]$ realizations. 
To achieve \eqref{eq:condition} for each
realization, start with an initial partition $\Delta{t}[1]$ and
then specify iteratively a new partition $\Delta{t}[k+1]$, from
$\Delta t [k]$, using the following refinement strategy:
\begin{equation}\label{eq:divide}
\begin{split}
&\mbox{{\bf for} each realization in the $j$-th batch  }\\
&\,\mbox{{\bf for} each time step }
n=1,\dots ,N_{\ssy A}[k]\\
&\quad\mbox{\bf if}~  r_n[k] \geq s_1 \tfrac{\tol_T}{\barN_{\ssy A}[j-1]},
\,\mbox{{\bf then}  divide $\Delta t_n[k]$ into $2$ uniform substeps.}\\
&\quad\mbox{{\bf else} let the new step be the same as the old} \\
& \quad\mbox{\bf endif} \\
&\,\mbox{\bf endfor}.\\
&\mbox{\bf endfor}.\\
\end{split}
\end{equation}

The refinement strategy \eqref{eq:divide} motivates the following
stopping criteria: for each realization of the $j$-th batch
\begin{equation}\label{eq:stop-cond1}
\mbox{\bf if}~ \left( \max_{1\leq n\leq
N_{\ssy A}[k]} r_n[k] < S_1 \tfrac{\tol_T}{\barN_{\ssy A}[j-1]}\right)
\quad\mbox{{\bf then}\quad stop,}
\end{equation}
where $S_1>\frac{2}{c}\,s_1$ with $c\approx\frac{1}{2}$,
see \cite{MSTZ}.
%
%
%
\par
Now we are ready for the detailed definition of the adaptive
algorithm with stochastic steps:
\medskip

{\tt Algorithm S}

\textbf{Initialization} Choose:
   \begin{enumerate}
   \item an error tolerance, $\tol \equiv \tol_S + \tol_T$,
   \item a number $N[1]$ of initial uniform steps $\Delta t[1]$ for $[0,T]$
    and set $\barN_{\ssy A} = N[1]$,
   \item a number $M[1]$ of initial realizations,
   \item a number $s_1 = 2$ in \eqref{eq:divide} and
           $c=\tfrac{1}{2}+\tfrac{1}{20}$ to compute $S_1$
           using $S_1>\frac{2}{c}\,s_1$, and
  \item a constant $c_0 \geq 1.65$ and an integer MCH$\ge 2$ to
    determine the number of realizations in \eqref{eq:choose_M}.
   \end{enumerate}

Set the iteration counter for batches $m=1$ and the stochastic
error $\ErrS[m] = + \infty$.

\textbf{Do while} ( $\ErrS[m] > \tol_S$ )

\hspace{0.5cm}\textbf{For} realizations $j = 1, \dots, M[m]$

\hspace{1cm}Set $k=1$ and $r[k]= + \infty$.

\hspace{1cm}Generate the jump times and their marks
$(\tau,Z)=\{(\tau_{\ell},Z_{\ell})\}_{\ell=1}^{{\widehat N}}$.
\par
\hspace{1cm}Start with the initial partition $\Delta t[k]$ and
    generate $\Delta W[k]$.
\par
\hspace{1cm}Compute $g(\barX(T))[J]$ and $N[J]$ by
{\tt routine Control-Time-Error}.

\hspace{0.5cm}\textbf{end-for}

\hspace{0.5cm}Compute the sample average
    $Eg \equiv \ma\left(g(\barX(T)); M[m]\right)$, the sample

\hspace{0.5cm}standard deviation
${\mathcal{S}}[m]\equiv{\mathcal{S}}(g(\barX(T)); M[m])$ and the a
posteriori bound

\hspace{0.5cm}for the statistical error
$\ErrS[m]\equiv\ErrS(g(\barX(T)), M[m])$ in \eqref{eq:static_bdd}.

\hspace{0.5cm}\textbf{if} ( $\mathtt{E}_S[m] > \tol_S$ )

\hspace{1cm}Compute $M[m+1]$ by {\tt change$\underline{\kenn}M$}%
        ($M[m]$, ${\mathcal{S}}[m]$, $\tol_S$; $M[m+1]$), cf.

\hspace{1cm}\eqref{eq:choose_M}, and update $\barN_{\ssy A} =\ma\left(N_{\ssy A}[J]; M[m] \right)$,
    where the random variable

\hspace{1cm}$N_{\ssy A}[J]$ is the final number of time steps
    on each realization.

\hspace{0.5cm}\textbf{end-if}

\hspace{0.5cm}Increase $m$ by 1.

\textbf{end-do}

Accept $Eg$ as an approximation of $E[g({{X}}(T))]$, since the
estimate of the

computational error is bounded by $\tol$.

\medskip

{\tt routine Control-Time-Error}$(\Delta t[k], \Delta W[k], r[k],
(\tau,Z); g(\barX(T))[J],N[J])$

\hspace{0.5cm}\textbf{Do while} ( $r[k]$ violates the stopping
    \eqref{eq:stop-cond1} )

\hspace{1.0cm}Compute the Euler approximation $\barX[k]$ in
Section~\ref{Kapote} and the error indicator

\hspace{1.0cm}$r[k]$ in \eqref{eq:indicator} using the error
density \eqref{eq:density} on $\Delta t[k]$ with the known

\hspace{1.0cm}Wiener increments $\Delta W[k].$

\hspace{1.0cm}\textbf{If} ( $r[k]$ violates the stopping
    \eqref{eq:stop-cond1} )

\hspace{1.5cm}\textbf{For} time steps $i = 1, \dots, N[k]$

\hspace{2.0cm}Do the refinement process \eqref{eq:divide} to
compute $\Delta t[k+1]$ from $\Delta t[k]$

\hspace{2.0cm} and compute $\Delta W[k+1]$ from $\Delta W[k]$
using Brownian bridges.

\hspace{1.5cm}\textbf{end-for}

\hspace{1.0cm}\textbf{end-if}

\hspace{1.0cm}Increase $k$ by 1.

\hspace{0.5cm}\textbf{end-do}

\hspace{0.5cm}Set the number of the final level $J = k-1$.

{\tt end of Control-Time-Error}
\section{Numerical Experiments}\label{sec:numerical_experiments}
\par
This sections shows numerical results from the implementation of
the a posteriori error approximation formula presented in
Section~\ref{sec:det_tstep} and of the adaptive algorithms
described in Section~\ref{sec:algo}. The programs we wrote uses
double-precision FORTRAN 77 and is based on the code written for
the numerical experiments in \cite{STZ}. For the numerical
simulation of the uniform distribution ${\mathcal U}(0,1)$ and the
normal distribution ${\mathcal N}(0,1)$, it applies a
double-precision modification of the functions {\tt ran1} and {\tt
gasdev} proposed in \cite{PTVF}, provided an initial seed which
must be a negative integer. In particular we use {\tt iseed} for
the simulation of the Wiener process increments, {\tt zseed} for
the simulation of the jump marks and {\tt tseed} for the
simulation of the jump times.
\par
To perform our computations, we consider a system of stochastic
differential equations of the form \eqref{eq:1.1} with:
$d=2$, $\ell_0=\ell_1=1$,
\begin{equation*}\label{chinese_exampleIII}
\gathered
a(t,x)=\big(\,-x_2,\,x_1+\tfrac{1}{2}\,\lambda(t)\,x_2\,\big),
\quad
b^1(t,x)=\Big(\,\sqrt{\tfrac{\lambda(t)}{1+t}}\,\sin(x_1),\,0\,\Big),\\
c(t,x,z)=\big(\,0,\,z\,\tfrac{\cos(x_1)}{\sqrt{t+1}}-x_2\,\big)
\endgathered
\end{equation*}
and a time dependent intensity $\lambda(t)=(1+t)^{-1}$.
The distribution for the jump marks is time dependent and such
that $E[Z_k^2] =1$. 
In particular, we use
\begin{equation*}
Z_k=\cos(2\pi\tau_k)+\sin(2\pi\tau_k)\,2\,\sqrt{3}\,(U_k-\tfrac{1}{2})
\end{equation*}
where $\{\tau_k\}$ are the jumps constructed in Remark \ref{Remarkk_1}
and $\{U_k\}$ is a sequence a sequence of ${\mathcal U}(0,1)$ i.i.d. random variables.
In this case, due to assumed for of $\lambda(t)$
the inverse function $\Lambda^{-1}$ is given explicitly by
$\Lambda^{-1}(s)=\exp(s)-1$.
This example is a generalization of Example 5.1 in \cite{LiuLi},
as here we admit a time dependent intensity of the underlying
Poisson process and also time dependent distribution of the marks.
Taking $g(x)=|x|^2$, $T=1$ and $X(0)=(0,0)$ the exact solution of the corresponding
weak approximation problem is given by the formula
\begin{equation*}
{E}\big[g(X(T))\big]=\left|X(0)\right|^2 +\int_0^T
\tfrac{\lambda(s)}{1+s}\;ds=1-\tfrac{1}{1+T} = \frac{1}{2}.
\end{equation*}
%
The value of the parameters needed in the simulations are
$\text{\rm MCH}=10$, $c_0=1.65$, $\text{\tt iseed}=-7$,
$\text{\tt zseed}=-101$ and $\text{\tt tseed}=-20$.
We note that the expected value of the number of jumps points
equals $\Lambda(T)=\Lambda(1)
=\text{\rm log}(2)\simeq
0.693$. Therefore, the computational cost of including the 
jump times of the process $X$ in our discretization is fairly low, see also the sampled values of $\max {\widehat N}$ in Table 5.3,
and it is asymptotically negligible as the required accuracy, $\tol$, tends to zero. 
%
%
\subsection{Deterministic time step algorithm}
First, we perform overkilling runs in order to test how realistic
is the a posteriori error approximation of the time discretization
error described in Theorem~\ref{thm:2.2}. The results, shown in
Table~5.1, show that the ratio of the computational error and its
computable approximation tends to 1 as the number of uniform time steps
$N$ increases. For each value of $N$, we choose the number of
realizations $M$ large enough in order to keep the total
statistical error at the level of 1\% of the size of the obtained
approximation of the time discretization error.
%
\input{table_overk}

\par
Table 5.2 contains the results of the {\tt Algorithm D} with  an
adaptive choice of both the deterministic time steps and the number of realizations.
The program starts with $M=100$ and $N=5$ subintervals as an
initial uniform partition of the time interval $[0,T]=[0,1]$. The
tolerance $\text{\rm TOL}=0.02$ is divided into
$\text{\rm TOL}_{\ssy S}=0.01333$, $\text{\rm TOL}_{\ssy T}=0.00444$ and
$\text{\rm TOL}_{\ssy T\!S}=0.00222$. When the algorithm stops, the
size of the total approximation error is less than
$2 \text{\rm TOL}$ according to the stopping
criterion \eqref{eq:stop-cond1_det} and it agrees with the
size of the computational error.
\newpage
\input{adpt1}

\subsection{Stochastic time step algorithm}
To observe the performance of the stochastic time steps {\tt
Algo\-rithm S}, we apply it for different values of
$\text{\rm TOL}$, starting with a number of realizations $M=100$ and a number of uniform time steps $N=5$. Table 5.3
contains the obtained results which show that {\tt Algorithm S}
is also effective in giving us an approximation of the quantity of
interest within the margin of $2\text{\rm TOL}$ 
due to the criterion \eqref{eq:stop-cond1_det}.
%
%
\begin{table*}[htdp]
\begin{center}
\begin{tabular}{|c|c|c|c|c|c|c|c|c|c|}
\hline
$\tol$  & $M$  &${\Cal A}(N_{\ssy A};M)$ & $\min N_{\ssy A}$
& $\max N_{\ssy A}$&${\Cal S}(N_{\ssy A};M)$& $\max {\widehat N}$
&    $E_{\ssy S}$ & ${\Cal E}_{\ssy C}$\\
\hline
$0.040$ &$3.2\!\times\!10^{3}$           &     8.3 &            5& 41&          5.2&  6
&   $2.5\!\times\!10^{-2}$&   $-2.2\!\times\!10^{-3}$\\
\hline
$0.020$ &$12.0\!\times\!10^{3}$&   10.9&   5&  73&          9.2&       6
&   $1.3\!\times\!10^{-2}$&   $-6.0\!\times\!10^{-3}$\\
\hline
$0.010$ &$47.9\!\times\!10^{3}$&   20.5&   10& 146&         16.9&  6
& $6.6\!\times\!10^{-3}$& $-1.0\!\times\!10^{-2}$\\
\hline
$0.005$ &$186.4\!\times\!10^{3}$ &   30.1&   10& 193&    32.2&   7
& $3.3\!\times\!10^{-3}$&   $-3.2\!\times\!10^{-3}$\\
\hline
\end{tabular}
\end{center}
\end{table*}

\par\noindent
\centerline{{\aut TABLE} 5.3. Adaptive choice of $M$ and $\Delta{t}$ with {\tt Algorithm S}.}
%
%
%
%
\par\medskip
\par\noindent
{\bf Acknowledgements.}
The work has been supported by: (i) the Swedish Research Council
for Enginee\-ring Science (TFR)  Grant\# 222-148, (ii) the Swedish
National Network in  Applied Mathematics (NTM), (iii) the Project
10.101 FCE - DINACYT, Uruguay, (iv) CSIC - Udelar, Human Resources
Program, and (v) the Facultad de Ciencias, Udelar, Montevideo,
Uruguay.
%
%
%
%
%
\appendix
\section{Discrete Dual Equations}\label{app:dualeqs}
This appendix section is dedicated to the determination
of the discrete dual functions $\varphi(t)\in\rset^d$,
$\varphi'(t)\in\rset^{d\times d}$ and
$\varphi''(t)\in \rset^{d\times d\times d}$
(see Theorem~\ref{thm:2.2} and
Theorem~\ref{thm:stoch_err_exp}),
where $t$ is a node of the (stochastic) partition
of the time interval $[0,T]$ which is used by
the Euler method (see Section~\eqref{Kapote}).
First, introduce the auxiliary functions ${\widehat A}_i$
and ${\widehat c}_i$, defined by
\begin{equation*}
\gathered
{\widehat A}_i(t_n,x)\equiv
x_i+\Delta t_n a_i(t_n,x)
+\Delta W_n^\ell b_i^\ell(t_n,x)\quad\forall\,x\in\rset^d,
\,\,i=1,\dots,d,\\
{\widehat c}_i(t,x,z)\equiv x_i+ c_i(t, x,z)
\quad\forall\,x\in\rset^d,\,\,\forall\,z\in{\Z},
\,\,i=1,\dots,d.\\
\endgathered
\end{equation*}
Then, for each realization, $\varphi$, $\varphi'$
and $\varphi''$ are constructed by the following
algorithm:
\par
\fbox{
\begin{minipage}[h]{14cm}
{\tt Dual backward time stepping algorithm.}
\par
\hspace{0.5cm} Set the initial backward values
\begin{equation}\label{eq:phiend}
\gathered
\varphi_i(t_{N_{\ssy A}})=\partial_i g\big(\barX(t_{N_{\ssy A}})\big),
\quad\varphi'_{ik}(t_{N_{\ssy A}})
=\partial_{ik} g\big(\barX(t_{N_{\ssy A}})\big),\\
\varphi''_{ikm}(t_{N_{\ssy A}})=\partial_{ikm} g\big(\barX(t_{N_{\ssy A}})\big).\\
\endgathered
\end{equation}
\hspace{0.5cm}\textbf{for} $n=N_{\ssy A}-1,\dots,0$

\hspace{1.0cm}\textbf{if} ($t_{n+1}$ is a jump time)
\textbf{then} set
$\Psistar=\big(t_{n+1}^{-},\barX(t_{n+1}^{-}),Z_{n+1}\big)$
and
\begin{equation}\label{eq:phijump}
\begin{split}
\varphi_i(t_{n+1}^{-})&=\partial_i\widehat{c}_j(\Psistar)\,\varphi_j(t_{n+1}),\\
\varphi'_{ik}(t_{n+1}^{-})&=\partial_i{\widehat c}_j(\Psistar)
\,\partial_k\widehat{c}_p(\Psistar) \,\varphi'_{jp}(t_{n+1})
+\partial_{ik}\widehat{c}_j(\Psistar)
\,\varphi_j(t_{n+1}),\\
\varphi''_{ikm}(t_{n+1}^{-})&=\partial_i\widehat{c}_j(\Psistar)\,
\partial_k\widehat{c}_p(\Psistar)\,
\partial_m\widehat{c}_r(\Psistar)\,
\varphi''_{jpr}(t_{n+1})\\
&\quad+\partial_{im}\widehat{c}_j(\Psistar)
\,\partial_k\widehat{c}_p(\Psistar)\,
\varphi'_{jp}(t_{n+1})\\
&\quad +\partial_i\widehat{c}_j(\Psistar)
\,\partial_{km}\widehat{c}_p(\Psistar)\,
\varphi'_{jp}(t_{n+1})\\
&\quad+\partial_{ik}\widehat{c}_j(\Psistar)
\,\partial_m\widehat{c}_p(\Psistar)\,
\varphi'_{jp}(t_{n+1})\\
&\quad
+\partial_{ikm}{\widehat c}_j(\Psistar)\,
\varphi_j(t_{n+1}),\\
\end{split}
\end{equation}
\hspace{1.5cm}\textbf{else} set
\begin{equation}\label{eq:phinojump}
\begin{split}
\varphi_i(t_{n+1}^{-})=\varphi_i(t_{n+1}),\quad
\varphi'_{ij}(t_{n+1}^{-})=\varphi'_{ij}(t_{n+1}),\quad
\varphi''_{ikm}(t_{n+1}^{-})=\varphi'_{ikm}(t_{n+1}).
\end{split}
\end{equation}
\hspace{1.0cm}\textbf{end-if}

\hspace{1.0cm} Set $\Rstar=(t_n,\barX(t_n))$ and
\begin{equation}\label{eq:phiinner}
\begin{split}
\varphi_i(t_n)&=\partial_i\widehat{A}_j(\Rstar)
\,\varphi_j(t_{n+1}^{-}),\\
\varphi'_{ik}(t_n)&=\partial_i{\widehat A}_j(\Rstar)
\,\partial_k{\widehat A}_p(\Rstar) \,\varphi'_{jp}(t_{n+1}^{-})
+\partial_{ik}{\widehat A}_j(\Rstar)
\,\varphi_j(t_{n+1}^{-}),\\
\varphi''_{ikm}(t_n)&=\partial_i\widehat{A}_j(\Rstar)\,
\partial_k\widehat{A}_p(\Rstar)\,
\partial_m\widehat{A}_r(\Rstar)\,
\varphi''_{jpr}(t_{n+1}^{-})\\
&\quad +\partial_{im}\widehat{A}_j(\Rstar)
\,\partial_k\widehat{A}_p(\Rstar)\,
\varphi'_{jp}(t_{n+1}^{-})\\
&\quad +\partial_i\widehat{A}_j(\Rstar)
\,\partial_{km}\widehat{A}_p(\Rstar)\,
\varphi'_{jp}(t_{n+1}^{-})\\
&\quad
+\partial_{ik}{\widehat A}_j(\Rstar)
\,\partial_m{\widehat A}_p(\Rstar)
\,\varphi'_{jp}(t_{n+1}^{-})
+\partial_{ikm}{\widehat A}_j(\Rstar)
\,\varphi_j(t_{n+1}^{-}).
\end{split}
\end{equation}
\hspace{0.5cm}\textbf{end-for}
\end{minipage}}
\par\noindent\smallskip\par
Observe that the estimate in Theorem~\ref{thm:2.2}
needs only the computation of $\varphi$ and $\varphi'$ and
it is therefore less expensive per realization.
Also, with respect to the actual implementation of the dual
backward time stepping it is useful to notice that the
blocks \eqref{eq:phijump} and \eqref{eq:phiinner}
differ only in the function call $\widehat{c}$ and
$\widehat{A}$.
%
%

%
\end{document}

%% file: table_overk.tex
%
%
%
%
\font\faan=eurm8
\def\Errr{\mathop{\text{\faan\char'105}}}
\def\EET{E_{\ssy T}}
\def\EES{E_{\ssy S}}
\def\EETS{E_{\ssy TS}}
\def\Cal{\mathcal}
%
%
%
%
%
\par\medskip\par\noindent
$$\vbox{\offinterlineskip
\hrule\halign{&\vrule#&
         \strut\hskip0.800truecm\hfil #\hfil\hskip0.800truecm &\vrule#\cr
%
%
%
}
\halign{&\vrule#&
         \strut\hskip2.540truecm\hfil #\hfil\hskip2.540truecm &\vrule#\cr
}
\halign{ &\vrule#&
        \strut\quad\hfil #\hfil\quad &\vrule#&
         \strut\quad\hfil #\hfil\quad &\vrule#&
          \strut\quad\hfil #\hfil\quad &\vrule#&
           \strut\quad\hfil #\hfil\quad &\vrule#&
            \strut\quad\hfil #\hfil\quad &\vrule#&
             \strut\quad\hfil #\hfil\quad &\vrule#&
              \strut\quad\hfil #\hfil\quad &\vrule#&
               \strut\quad\hfil #\hfil\quad &\vrule#\cr
height4pt&\omit&&\omit&&\omit&&\omit&&\omit&\cr
&   $\boxed{\text{\tt Alg. D}}$
&&  $c_0=1.65$
&&  $\text{\tt iseed}=-7$
&&  $A:=\tfrac{\EET-\EES-\EETS}{{\Cal E}_c}$
&&  $B:=\tfrac{\EET+\EES+\EETS}{{\Cal E}_c}$&\cr
height4pt&\omit&&\omit&&\omit&&\omit&&\omit&\cr
\noalign{\hrule}
height4pt&\omit&&\omit&&\omit&&\omit&&\omit&\cr
&  $N$
&& $M$
&& $\EET$
&& $\EES+\EETS$
&& $[\,\min\{A,B\},\ \ \max\{A,B\}\,]$ &\cr
height4pt&\omit&&\omit&&\omit&&\omit&&\omit &\cr
\noalign{\hrule}
height3pt&\omit&&\omit&&\omit&&\omit&&\omit &\cr
&  $5$
&& $10\times10^6$
&& $-0.0602$
&& $5.87\times 10^{-4}$
&& $[1.026,\ \ 1.046]$ &\cr
height3pt&\omit&&\omit&&\omit&&\omit&&\omit &\cr
\noalign{\hrule}
height3pt&\omit&&\omit&&\omit&&\omit&&\omit &\cr
&  $10$
&& $50\times10^6$
&& $-0.0314$
&& $2.33\times 10^{-4}$
&& $[1.019,\ \ 1.035]$ &\cr
height3pt&\omit&&\omit&&\omit&&\omit&&\omit &\cr
\noalign{\hrule}
height3pt&\omit&&\omit&&\omit&&\omit&&\omit &\cr
&  $20$
&& $100\times 10^6$
&& $-0.0159$
&& $1.54\times 10^{-4}$
&& $[1.008,\ \ 1.028]$  &\cr
height3pt&\omit&&\omit&&\omit&&\omit&&\omit &\cr
}
\hrule}$$
\par\noindent
\centerline{{\aut TABLE} 5.1. Computing the  efficient index
of the {\tt Algorithm D}.}
\par\medskip\par

%% file: adpt1.tex
%
%
%
%
%
%
%
%
%
%
%
%
$$\vbox{\offinterlineskip \hrule
\halign{&\vrule#&
         \strut\hskip1.145truecm\hfil #\hfil\hskip1.145truecm &\vrule#\cr
}

  \halign{&\vrule#&
         \strut\quad\hfil #\hfil\quad &\vrule#&
          \strut\quad\hfil #\hfil\quad &\vrule#&
           \strut\quad\hfil #\hfil\quad &\vrule#&
            \strut\quad\hfil #\hfil\quad &\vrule#&
             \strut\quad\hfil #\hfil\quad &\vrule#&
              \strut\quad\hfil #\hfil\quad &\vrule#&
               \strut\quad\hfil #\hfil\quad &\vrule#&
                \strut\quad\hfil #\hfil\quad &\vrule#&
                 \strut\quad\hfil #\hfil\quad &\vrule#\cr
height5pt&\omit&&\omit&&\omit&&\omit&&\omit&&\omit&&\omit&\cr
& \boxed{\text{\tt Alg. D}}
&& $ $
&& $c_{_0}=1.65$
&& $\text{\rm MCH}=10$       && $\text{\tt iseed=-7}$
&& $ \text{\rm TOL}=0.02$
&& $ $ &\cr
height5pt&\omit&&\omit&&\omit&&\omit&&\omit&&\omit&&\omit&\cr
\noalign{\hrule}
height5pt&\omit&&\omit&&\omit&&\omit&&\omit&&\omit&&\omit&\cr
& Iter. && $N$ && $M$ && ${\Cal E}_c$ && $\EET$&&$\EETS$
&& $\EES$ &\cr
height5pt&\omit&&\omit&&\omit&&\omit&&\omit&&\omit&&\omit&\cr
\noalign{\hrule}
height3pt&\omit&&\omit&&\omit&&\omit&&\omit&&\omit&&\omit&\cr
   & 1  &&  5 &&  100 && -0.03272 &&  -0.05889 && 0.02112
   && 0.14602 &\cr
height3pt&\omit&&\omit&&\omit&&\omit&&\omit&&\omit&&\omit&\cr
\noalign{\hrule}
height3pt&\omit&&\omit&&\omit&&\omit&&\omit&&\omit&&\omit&\cr
   & 2  &&  5 &&  1000  && -0.04844 &&  -0.06005 && 0.00684
   && 0.04835  &\cr
height3pt&\omit&&\omit&&\omit&&\omit&&\omit&&\omit&&\omit&\cr
\noalign{\hrule}
height3pt&\omit&&\omit&&\omit&&\omit&&\omit&&\omit&&\omit&\cr
   & 3  && 5 &&  10000 && -0.06592 &&  -0.06070 && 0.00232
   && 0.01675 &\cr
height3pt&\omit&&\omit&&\omit&&\omit&&\omit&&\omit&&\omit&\cr
\noalign{\hrule}
height3pt&\omit&&\omit&&\omit&&\omit&&\omit&&\omit&&\omit&\cr
   & 4  &&  5 && 12088  && -0.05215 &&  -0.05925 && 0.00204
   && 0.01459  &\cr
height3pt&\omit&&\omit&&\omit&&\omit&&\omit&&\omit&&\omit&\cr
\noalign{\hrule}
height3pt&\omit&&\omit&&\omit&&\omit&&\omit&&\omit&&\omit&\cr
   & 5  && 10 && 12088 && -0.03738 &&  -0.03196 && 0.00107
   && 0.01403 &\cr
height3pt&\omit&&\omit&&\omit&&\omit&&\omit&&\omit&&\omit&\cr
\noalign{\hrule}
height3pt&\omit&&\omit&&\omit&&\omit&&\omit&&\omit&&\omit&\cr
   & 6  && 20 && 12088 && -0.02585 &&  -0.01633 && 0.00054
   && 0.01369 &\cr
height3pt&\omit&&\omit&&\omit&&\omit&&\omit&&\omit&&\omit&\cr
\noalign{\hrule}
height3pt&\omit&&\omit&&\omit&&\omit&&\omit&&\omit&&\omit&\cr
   & 7  && 20 && 14122  && -0.02559 &&  $-$ && $-$
   && 0.01261 &\cr
height3pt&\omit&&\omit&&\omit&&\omit&&\omit&&\omit&&\omit&\cr
%
 }\hrule}$$
\par\noindent
\centerline{{\aut TABLE} 5.2. Adaptive choice of $M$
and $\Delta{t}$ with {\tt Algorithm D}.}
\par\medskip\par